\documentclass{amsart}
\usepackage[dvips]{graphics} 
\input diagrams
\diagramstyle[PostScript=dvips]

\newcommand{\stack}[1]{{\mathfrak #1}}

\newcommand{\mg}{\stack{M}_g}
\newcommand{\mgbar}{\overline{\stack{M}_g}}

\newcommand{\sgr}{{\stack S}_{g}^{1/r}}

\newcommand{\sgrnm}{{\stack S}^{1/r,\mathbf m}_{g,n}}

\newcommand{\sgrbar}{\overline{\stack S}_{g}^{1/r}}

\newcommand{\sgrnmbar}{\overline{\stack S}^{1/r,\mathbf m}_{g,n}}

\newcommand{\stacksgr}{{\stack S}_{g}^{1/r}}
\newcommand{\stacksgrnmbar}{\overline{\stack S}^{1/r,\mathbf m}_{g,n}}
\newcommand{\stacksgrbar}{\overline{\stack S}_{g}^{1/r}}
\newcommand{\stacksgsbar}{\overline{\stack S}_{g,s}}
\newcommand{\stacksgrnm}{{\stack S}^{1/r,\mathbf m}_{g,n}}
\newcommand{\stackmg}{\stack{M}_g}
\newcommand{\stackmgbar}{\overline{\stack{M}}_g}

\newcommand{\stackmgnbar}{\overline{\stack{M}}_{g,n}}
\newcommand{\stacksgtwo}{{\stack S}_{g,2}} 
\newcommand{\spec}{\text{Spec\,}}
\newcommand{\proj}{\text{Proj\,}}
\newcommand{\pic}{\text{Pic\,}}
\newcommand{\jac}{\text{Jac\,}}
\newcommand{\tensor}{\otimes}
\newcommand{\cross}{\times}
\newcommand{\irightarrow}{\rTo^{\sim}}
\theoremstyle{plain}
\newtheorem{theorem}{Theorem}[subsection]
\newtheorem{proposition}[theorem]{Proposition}
\newtheorem{lemma}[theorem]{Lemma}
\newtheorem{corollary}[theorem]{Corollary}

\newtheorem{conjecture}[theorem]{Conjecture}
\theoremstyle{remark}

\newtheorem{nb}[theorem]{Note}

\theoremstyle{definition}
\newtheorem{defn}[theorem]{Definition}

\theoremstyle{plain}
\newcommand{\cx}{{\mathcal X}}

\newcommand{\ck}{{\mathcal K}}
\newcommand{\cl}{{\mathcal L}}
\newcommand{\ce}{{\mathcal E}}
\newcommand{\cf}{{\mathcal F}}
\newcommand{\co}{{\mathcal O}}

\newcommand{\ch}{{\mathcal H}}

\newcommand{\cy}{{\mathcal Y}}
\newcommand{\cm}{{\mathcal M}}
\newcommand{\cw}{{\mathcal W}}
\newcommand{\cn}{{\mathcal N}}

\newcommand{\sheafhom}{\ch\kern-.15em om}

\newcommand{\bm}{\mathbf{m}}          %

\newcommand{\me}{\mathfrak{E}_{r}}
\newcommand{\te}{\tilde{\ce}_{r}}

\begin{document}
\setcounter{tocdepth}{2}

\title[Picard group of spin curve moduli]{The Picard group of the moduli of
higher spin curves} 

\author{Tyler J. Jarvis}
\address{Department of Mathematics\\ Brigham Young University\\Provo, UT 84602}
\email{jarvis@math.byu.edu}
\thanks{This material is based on work supported in part by the National
Science Foundation under Grant No.~DMS-9501617 and by the 
National Security Agency under Grant No.~MDA904-99-1-0039}
\date{\today}

\begin{abstract}
  
  This article treats the Picard group of the moduli (stack)
  $\sgr$ of $r$-spin curves and its compactification $\sgrbar$.  Generalized
  spin curves, or $r$-spin curves are a natural generalization
  of $2$-spin curves (algebraic curves with a theta-characteristic),
  and have been of interest lately because they are the subject of a
  remarkable conjecture of E. Witten, and because of the 
  similarities between the intersection theory of these moduli 
  spaces and that of the moduli of stable maps.
  
 We generalize results of Cornalba, describing and giving
  relations between many of the elements of the Picard group of the stacks $\stacksgr$
  and $\stacksgrbar$. These relations are important in the proof of the genus-zero case of 
  Witten's conjecture given in \cite{jkv}. We use these relations to show that when $2$ or $3$ divides $r$, then
  $\pic \sgr$ has non-zero torsion.  And finally, we work out 
  some specific examples.

\end{abstract}

\maketitle

\section{Introduction}

In this article we study the Picard group of the moduli  (or 
rather the stack) $\sgrbar$ of higher spin curves, or $r$-spin 
curves, over the Deligne-Mumford stack of stable curves 
$\stackmgbar$.  Smooth $r$-spin curves consist of a smooth 
algebraic curve $X$, a line bundle (invertible coherent sheaf) 
$\cl$, and an isomorphism from the $r$th tensor power 
$\cl^{\tensor 
  r}$ to the cotangent bundle $\omega_X$.  The
compactification of the stack of smooth spin curves uses {\em 
stable spin curves}.  These consist of a stable curve $X$, and 
for any divisor $d$ of $r$, a rank-one, torsion-free sheaf 
$\ce_d$ which is almost a $d$th root of the canonical (relative 
dualizing) sheaf.  This is made precise in 
Definition~\ref{stablespin}. 

The stack of spin curves  provides  a finite cover of the stack 
of stable curves.  Although this stack has many similarities to 
the stack of  stable maps, including the existence of classes 
analogous to Gromov-Witten classes, and an associated 
cohomological field theory \cite{jkv}, it is not the stack of 
stable maps into any variety \cite[\S5.1]{jkv}.  

These moduli spaces are especially interesting because of a 
conjecture of E. Witten relating the intersection theory on the 
moduli space of $r$-spin curves and KdV (Gelfand-Dikii) 
hierarchies of order $r$ 
\cite{witten:N-matrix-model,witten:r-spin-conj}. This conjecture 
is a generalization of an earlier conjecture of his, which was 
proved by Kontsevich (see 
\cite{kontsevich:thesis,kontsevich:matrix-airy} and 
\cite{looijenga:kontsevich}).  As in the case of Gromov-Witten 
theory, one can construct a cohomology class $c^{1/r}$ and a 
potential function from the intersection numbers of $c^{1/r}$ and 
the tautological $\psi$ classes associated to the universal 
curve. Witten conjectures that the potential of the theory 
corresponds to the tau-function of the order-$r$ Gelfand-Dikii 
($KdV_r$) hierarchy.   In genus zero, the Witten conjecture is 
true \cite{jkv}, and the relations of Theorem~\ref{mainthm} play 
a role in the proof. 

Construction of the stack $\sgr$, and its compactification 
$\sgrbar$, was done in \cite{cornalba:theta} for $r=2$ and in 
\cite{jarvis:geometry} for all $r\ge 2$.  This article focuses on 
giving a description of the Picard group of $\sgr$ and $\sgrbar$. 

 \subsection{Overview and Outline of the Paper}

 In Section~\ref{11} and \ref{12} we give definitions of
 smooth and stable spin curves and some examples.  The examples of
 Section~\ref{222} are typical of the spin curves that arise over the
 boundary divisors of $\pic \mgbar$. In Section~\ref{geometry} we recall the basic properties of the 
moduli spaces from \cite{jarvis:geometry}.  
 In Section~\ref{2} we treat the Picard group of these spaces.
 
 The Picard group that we will work with in this paper is the Picard
 group of the stack $\pic \stacksgrbar$ (also called $\pic_{fun}$), as defined in
 \cite{mumford:picard} or \cite{harris-mumford}.  The exact definition of
 the Picard group is given in Section~\ref{21}, and the definitions of
 the boundary divisors and the tautological elements of $\pic \stacksgrbar$
 are given in Section~\ref{div-basic}.  We show in
 Section~\ref{224} that the boundary divisors and the Hodge class are
 independent in $\pic \stacksgrbar$.  This turns out to be a useful step
 toward showing that torsion exists in $\pic \stacksgr$.
 
 In Section~\ref{23} we compute the main relations between
 elements of the Picard group, and we show some of the consequences of
 those relations.  The general results are given in
 Theorems~\ref{mainthm} and ~\ref{mainthm}.bis, and some special cases
 of those results are given in Corollary~\ref{table}.  In the special
 case that $r=2$, these results reduce to  those of
 \cite{cornalba:theta} and \cite{cornalba:spin-pic}.
 
 In section \ref{tors} we discuss the existence of torsion in $\pic \stacksgr$ when $2$ or $3$ divides $r$.  In
 particular, Proposition~\ref{torsion} shows that when $2$ divides
 $r$, $\pic \stacksgr$ has $4$-torsion, and when $3$ divides $r$, $\pic
 \stacksgr$ has $3$-torsion.   This is in stark 
contrast to the case of $\pic \mg$, and $\pic 
\mgbar$, which are known to be free for $g \geq 3$. 
\cite{arbarello-cornalba:mg-pic}.  

Finally, in Section \ref{elliptic} we work out examples for genus 
$1$ and general $r$, and for $r=2$ and general genus. 

 \subsection{Previous Work}
 
 The Picard group of the moduli space of curves is now fairly well
 understood, thanks primarily to the work of Arbarello-Cornalba
 \cite{arbarello-cornalba:mg-pic} and Harer \cite{harer:mg-h2}.  Most
 of the progress toward understanding the Picard group of the moduli
 of $2$-spin curves is also due to Cornalba \cite{cornalba:theta,
   cornalba:spin-pic} and Harer \cite{harer:spin-pic-rat}.
 Section~\ref{23} on relations between classes in the Picard group
 is strongly motivated by Cornalba's work in \cite{cornalba:theta} and
 \cite{cornalba:spin-pic}.
 
In \cite{jarvis:geometry} and \cite{jarvis:spin} several 
different compactifications  of the moduli of $r$-spin curves are 
constructed.  The best-behaved of these compactifications, and  
the one we will use here, is the stack of what are called {\em 
stable $r$-spin curves} or just {\em 
   $r$-spin curves} in \cite{jarvis:geometry}.  In the case of prime $r$, 
   these are the same as the \emph{pure spin curves} of 
   \cite{jarvis:spin}.

 \subsection{Conventions and Notation}
 
 By a {\em curve} we mean a reduced,
 complete, connected, one-dimensional scheme over an algebraically closed field $k$.  A {\em
   semi-stable curve of genus $g$} is a curve with only ordinary
 double points such that $H^1(X,\co_X)$ has dimension $g$.  And an
 {\em $n$-pointed stable curve} is a semi-stable curve $X$ together
 with an ordered $n$-tuple of non-singular points $(p_1,\dots,p_n)$,
 such that at least three marked points or double points of $X$ lie on
 every smooth irreducible component of genus $0$, and at least one
 marked point or double point of $X$ lies on every smooth component of
 genus one.  A family of stable (or semi-stable) curves is a flat,
 proper morphism $ \cx \rightarrow T $ whose geometric fibres $\cx_t$
 are (semi) stable curves.  Irreducible components of a semi-stable
 curve which have genus $0$ (i.e., are birational to $\mathbb P^1$)
 but which meet the curve in only two points will be called {\em
   exceptional} curves.
 
 By {\em line bundle} we mean an invertible (locally free of rank one)
 coherent sheaf.  By canonical sheaf we will mean the relative
 dualizing sheaf of a family of curves $f:\cx \rightarrow T$, and this
 sheaf will be denoted $\omega_{\cx/T}$ or $\omega_f$.  Note that for a
 semi-stable curve, the canonical sheaf is a line bundle.  When $T$ is
 the spectrum of an algebraically closed field, we will often write $X$ for $\cx$, and $\omega_X$ for
 $\omega_{X/T}$.

\section{Spin Curves} \label{spin-curves-sect}
        \subsection{Smooth Spin Curves}\label{11}

For any $g \geq 0$, fix a positive integer $r$ with the property 
that $r$ divides $2g-2$.  We define a smooth {\em $r$-spin curve} 
to be a triple $(X,\cl,b)$ of a smooth curve $X$ of genus $g$, a 
line bundle $\cl$,  and an isomorphism $b$, of the $r$th tensor 
power of $\cl$ to the canonical bundle $\omega_{X}$ of $X$; that 
is, $b:\cl^{\otimes r} \irightarrow 
\omega_{X}$.  For a given $X$, any $(\cl,b)$ making $(X,\cl,b)$ 
into a spin curve will be called an {\em $r$-spin structure}.  
Families of smooth $r$-spin curves are triples $(\cx/T, \cl,b)$ 
of a family $\cx/T$ of smooth curves, a line bundle $\cl$, and an 
isomorphism $b$, which induces an $r$-spin structure on each 
geometric fibre of $\cx/T$. 

{\bf Example: $\mathbf{r=2}$.}

A $2$-spin curve is what has classically been called a spin curve  
(a curve with a theta-characteristic $\cl$) \emph{with an 
explicit isomorphism} $b:\cl^{\otimes 2} \irightarrow \omega$.

For any choice of $\mathbf m = (m_1,m_2,m_3,\dots,m_n)$, such 
that $r$ divides $2g-2-\sum m_i$, we may also define an 
$n$-pointed $r$-spin curve of type $\mathbf{m}$  as a triple 
$((X,(p_1,p_2,\dots,p_n)), \cl,b)$ such that 
$(X,(p_1,p_2,\dots,p_n))$ is a smooth, $n$-pointed curve, and $b$  
is an isomorphism from $\cl^{\otimes r}$ to $ \omega_{X}(-\sum 
m_i p_i)$. Families of $n$-pointed $r$-spin curves are defined 
analogously. 

Two spin structures $(\cl,b)$ and $(\cl',b')$ are isomorphic   if 
there exists an isomorphism between the bundles $\cl$ and $\cl'$ 
which respects the homomorphism $b$ and $b'$.  Over a fixed curve 
$X$, any two $r$-spin structures $(\cl,b)$ and $(\cl,b')$ which 
differ only by their isomorphism $b$ or $b'$ must be isomorphic,  
provided the base is algebraically closed.   The set $\sgr[X]$ of 
isomorphism classes of $r$-spin structures on $X$ is a principal 
homogeneous space over the $r$-torsion $\jac_r X$ of the Jacobian 
of $X$; thus it has $r^{2g}$ elements in it. 

Similarly, two spin curves $(X,\cl,b)$ and $(X',\cl',b')$ are 
isomorphic if there is an isomorphism $\tau: X 
\irightarrow X'$ and an isomorphism of spin structures $i: \cl 
\irightarrow 
\tau^* 
\cl$. 

{\bf Example 2: $\mathbf{g=1, \, r\ge 2, \mathbf{m}=0}$.}

If $n \ge 1$ and $\mathbf m = \mathbf 0$, then, up to 
isomorphism, an $n$-pointed $r$-spin curve of genus $1$ is just 
an $n$-pointed curve of genus  $1$ with a point of order $r$ on 
the curve.  However, the automorphisms of the underlying curve 
identify some of these $r$-spin structures.  In particular, when 
$n=1$ and $r$ is odd, the elliptic involution acts freely on all 
the non-trivial $r$-spin structures, thus there are only 
$1+(r^2-1)/2$ isomorphism classes of $r$-spin curves over a 
generic $1$-pointed curve of genus $1$.

Let $\sgr$ denote the stack of smooth $r$-spin curves of genus 
$g$.  And let $\sgrnm$ denote the stack of $n$-pointed spin 
curves (for a given $\mathbf m$).  When $g=1$, and $n=1$, we   
will also write $\stack{S}_{1}^{1/r}$ to denote the stack 
$\stack{S}^{1/r, \mathbf{0}}_{1,1}$.  If $X$ has no 
automorphisms, the sets $\sgr[X]$ and $\sgrnm[X]$, as defined 
above, are just the fibres of $\sgr$, or $\sgrnm$ over the point 
corresponding to $X$ in $\mg$, or in $\stack{M}_{g,n}$, 
respectively.

\subsection{Stable Spin Curves}\label{12}
        
To compactify the moduli of spin curves, it is necessary to 
define a spin structure for stable curves.  To do this we need 
not just line bundles, but also rank-one torsion-free sheaves 
over stable curves.  Some additional structure, as given in 
Definitions~\ref{root} and \ref{systems}, is also necessary to 
ensure that the compactified moduli space $\sgrnmbar$ is 
separated and smooth. 

  \subsubsection{Definitions}

  To begin we need the definition of torsion-free sheaves.
 \begin{defn}
  A {\em relatively torsion-free sheaf} (or just torsion-free sheaf)
  on a family of stable or semi-stable curves $f: \cx \rightarrow T$
  is a coherent $\co_{\cx}$-module $\ce$ that is flat over $T$, such
  that on each fibre $\cx_t =\cx \times_T \spec k(t)$ the induced
  $\ce_t$ has no associated primes of height one.
 \end{defn}

We will only be concerned with rank-one torsion-free sheaves.  Such
sheaves are called {\em admissible} by Alexeev
\cite{alexeev:compact-jac} and {\em sheaves of pure dimension $1$} by Simpson
\cite{simpson:rep-fund-group}.  Of course, on the open set where $f$
is smooth, a torsion-free sheaf is locally free.

\begin{nb}
 It is well-known and easy to check that if a rank-one, 
torsion-free sheaf $\ce$ is not locally free (also called {\em 
singular}) at a node $p$ of $X$, then the completion 
$\hat{\co}_{X,p}$ of the local ring of $X$ near $p$ is isomorphic  
to $A = k[[x,y]]/xy$, and $\ce$ corresponds to an $A$-module $E 
\cong xk[[x]] \oplus yk[[y]]=<\zeta_1,\zeta_2|y\zeta_1=x\zeta_2=0>$ \cite[Prop. 11.3]{seshadri:fibres}. 
\end{nb}

\begin{defn}\label{root}\label{rootdef}
Given an $n$-pointed, semi-stable curve $(X, p_1, \dots, p_n)$, 
and a rank-one, torsion-free sheaf $\ck$ on $X$, and given an  
$n$-tuple $\mathbf{m} = (m_1, \dots, m_n)$ of integers, we denote 
by $\ck(\bm)$ the sheaf $\ck \otimes \co(-\sum m_i p_i)$. 

A \emph{$d$th root of $\ck$ of type $\mathbf{m}$} is a pair 
$(\ce, b)$ of a rank-one, torsion-free sheaf $\ce$, and an 
$\co_X$-module homomorphism $b: 
\ce^{\tensor d} 
 \rTo \ck (\bm)$ with the following properties:
\begin{enumerate}
\item $d \cdot \deg \ce = \deg \ck-\sum m_i$
\item $b$ is an isomorphism on the locus of $X$ where $\ce$ is 
locally free
\item \label{coker} for every point $p \in X$ where $\ce$ is not free, the 
length of the cokernel of $b$ at $p$ is $d-1$.
\end{enumerate} 
\end{defn}

Unfortunately, the moduli space of stable curves with $d$th roots 
of a fixed sheaf $\ck$ is not smooth when $d$ is not prime, and 
so we must consider not just roots of a bundle, but rather a 
coherent net of roots.  This additional structure suffices to 
make the stack of stable curves with coherent root nets smooth  
\cite{jarvis:geometry}. 

\begin{defn}\label{systems}
Given a semi-stable $n$-pointed curve $(X, p_1, \dots, p_n)$, and 
a rank-one, torsion-free sheaf $\ck$ on $X$, and an $n$-tuple 
$\mathbf{m}$, a 
\emph{coherent net of roots of type $\mathbf{m}$ for $\ck$} 
is a collection $\{\ce_d, c_{d,d'}\}$ consisting of a rank-one 
torsion-free sheaf $\ce_d$ for each $d$ dividing $r$, and a 
homomorphism $c_{d,d'}: \ce^{\otimes d/d'}_d \rTo \ce_{d'}$ for 
each $d'$ dividing $d$ with the following properties: 
\begin{itemize}
\item $\ce_1= \ck$, and $c_{d,d} =1$ for each $d$ dividing $r$.
\item For every pair of divisors $d'$ and $d$ of $r$ such that $d'$ divides $d$, let $\mathbf{m}'$ 
be the $n$-tuple $(m'_1, \dots, m'_n)$ such that $m'_i$ is the 
smallest, non-negative integer congruent to $m_i 
\operatorname{mod}(d/d')$.  The  $\co_X$-module homomorphism 
$c_{d,d'} 
: 
\ce^{\tensor d/d'}_{d} 
\rTo 
\ce_{d'}$, must make $\ce_d$ into a $d/d'$-root of $\ce_{d'}$, 
of type $\mathbf{m}'$, such that all these maps are compatible.  
That is, the diagram 
$$
\begin{diagram}
(\ce^{\tensor d/d'}_{d})^{\tensor d'/d''} & \rTo^{(c_{d,d'})^{\tensor d'/d''}} & \qquad\ce^{\tensor d'/d''}_{d'}\\ 
&  \rdTo^{c_{d,d''}}  & \dTo c_{d',d''} \\
& &  \ce_{d''}\\
\end{diagram}
$$
commutes for every $d''|d'|d|r$.  
\end{itemize}
\end{defn}

If $r$ is prime, then a coherent net of $r$th-roots is simply an 
$r$th root of $\ck$.  Moreover, if $\ce_d$ is locally free, then 
up to isomorphism $\ce_d$ uniquely determines all $\ce_{d'}$ and 
all $c_{d,d'}$ such that $d'|d$. 

\begin{defn}\label{stablespin}
A \emph{stable, $n$-pointed, $r$-spin curve of type $\bm = (m_1, 
\dots, m_n)$} is an $n$-pointed, stable curve $(X, 
p_1, 
\dots, p_n)$ and a coherent net of $r$th roots of $\omega_X $ of type $\mathbf{m}$, 
where $\omega_X$ is the canonical (dualizing) sheaf of $X$.  An 
$r$-spin curve is called \emph{smooth} if $X$ is smooth. 
\end{defn}

Note that this definition of a smooth $r$-spin curve differs from 
that of Section~\ref{11}, in that a spin curve carries the 
additional data of explicit isomorphisms  $\ce^{\otimes d/d'}_d 
\rTo \ce_{d'}$; however, for smooth curves,  and indeed, whenever 
$\ce_r$ is locally free, $\ce_r$ and $c_{r,1}$ completely 
determine the spin structure, up to isomorphism. 

\begin{defn}
An \emph{isomorphism of $r$-spin curves} from $(X, p_1, \dots, 
p_n, 
\{\ce_d, c_{d, d'}\})$ to $(X', p'_1, \dots, p'_n, \{\ce'_d, c'_{d, 
d'}\})$ is an isomorphism of pointed curves $$\tau : (X, p_1, 
\dots, p_n)  \rTo^~ (X', p'_1, \dots, p'_n)$$ and a system of 
isomorphisms $\{\beta_d : \tau^* \ce'_d  \rTo^~ \ce_d\},$ with 
$\beta_1$ the canonical isomorphism $\tau^* \omega_{X'}(-\sum_i 
m_i p'_i)  \rTo^~ 
\omega_X(-\sum m_ip_i),$ and such that the $\beta_d$ are compatible with 
all of the maps $c_{d,d'}$ and $\tau^*c'_{d,d'}$.
\end{defn}

The definition of families of $r$-spin curves is relatively 
technical and unenlightening.  For the details of those 
definitions see \cite{jarvis:geometry}.  For our purposes, it 
will suffice to know the basic properties of the stack of 
$r$-spin curves from \cite{jarvis:geometry} as given in Section 
\ref{geometry} 

\subsubsection{Alternate Description of Stable Spin Curves}

The following characterization of spin curves in terms of line bundles
on a partial normalization of the underlying curve is very useful and
helps illustrate the nature of stable spin curves.

Consider a stable spin curve $(X, \{\ce_d, c_{d,d'}\})$ and the 
partial normalization $\tilde{X} \rTo^{\pi} X$ of $X$ at each of 
the singularities of $\ce_r$ (i.e., the nodes of $X$ where 
$\ce_r$  fails to be locally free).  The completion 
$\hat{\co}_{X,p}$ of the local ring of $X$ near a singularity $p$ 
of $\ce_r$ is isomorphic to $A = k [[x,y]]/xy$, and $\ce_r$ 
corresponds to an  $A$-module $E \cong xk[[x]] \oplus y k[[y]]= 
<\zeta_1, \zeta_2| x 
\zeta_2=y\zeta_1=0>$.  The homomorphism $c_{r,1}:\ce^{\otimes r}_r \rTo \omega$ 
corresponds to a homomorphism of $A$-modules $E^{\otimes r} \rTo 
A$ of the form $\zeta^r_1 \mapsto x^u$, $\zeta^r_2 \mapsto y^v$, 
and $\zeta^i_1 \zeta^{r-i}_2 \mapsto 0$ for $0 <i<r$.  The 
condition on the cokernel (Definition \ref{rootdef}(\ref{coker})) 
implies that $u+v=r$.  The pair $\{u,v\}$ is called the 
\emph{order} of $c_{r,1}$ at $p$.  If $\ce_r$ is locally free at $p$, 
then $c_{r,1}$ is an isomorphism and the order is $\{0,0\}$.  

Let  $\tilde{A}= k[[x]] \oplus k [[y]] \cong 
\hat{\co}_{\tilde{X},p^+} \oplus \hat{\co}_{\tilde{X},p^-}$, where 
$\{p^+, p^-\}$ is the inverse image $\pi^{-1}(p)$ of the 
normalized node. The pullback $\pi^*\ce_r$ of $\ce_r$ corresponds 
to $E 
\otimes_A 
\tilde{A}$, which is no longer torsion-free; but if $\pi^{\natural} \ce_r:= \pi^* \ce_r/torsion$, 
which corresponds to a free $\tilde{A}$-module $\tilde{E}$, then
$c_{r,1}$ induces an isomorphism $\tilde{c}_{r,1}: 
\pi^{\natural} \ce^{\otimes r}_r \irightarrow \pi^* \omega_X (-u p^+-vp^-)=
\omega_{\tilde{X}}((1-u)p^++(1-v)p^-)$.  

Conversely, given  a partial normalization $\pi : \tilde X 
\rightarrow X$ 
 with $\pi^{-1}(q_i) = \{ q_i^{+}, q_i^{-}\}$, the inverse images of 
the singular points $q_i$, and given integers $u_i \in (0,r)$ and 
$v_i 
\in (0,r)$ for each normalized singularity $q_i$ such that 
$u_i + v_i = r$, consider an $r$th root $\cl$ of the line bundle 
$\pi^* \omega_{X} 
\tensor \co_{\tilde X}(-\sum( u_i q_i^{+} + v_i q_i^{-})) \cong
\omega_{\tilde X} \tensor \co_{\tilde X} (-\sum ( (u_i-1) q_i^{+} +
(v_i-1) q_i^{-})) $. Of course for such an $\cl$ to exist, $u_i$ 
and $v_i$ must be chosen to make the degree of the $r$th power of 
$\cl$ divisible by $r$.  From this $r$th root and partial 
normalization we can create an $r$th root of $\omega_X$ on $X$ by 
taking $\ce_r 
= 
\pi_*{\cl}$, and by taking $c_{r,1}$ to be the map induced by adjointness 
from the composite $\cl^{\tensor r} 
\irightarrow 
\pi^* \omega_{X} \tensor 
\co_{\tilde X}(-\sum u_i q_i^{+} + v_i q_i^{-}) \hookrightarrow \pi^*
\omega_{X}.$   Thus $r$th roots of $\omega_X$ are in one-to-one correspondence
with $r$th roots of $\omega_{\tilde{X}}(-(u-1)q^+_i-(v-1)q^-_i)$ 
on $\tilde{X}$.  

Moreover, if $\{u,v\}$ is the order of $c_{r,1}:\ce^{\otimes r}_r 
\rTo \omega_X$ at $p$, then for each $d$ dividing $r$, the order of 
$(\ce_d, c_{d,1})$ is $\{u_d,v_d\}$ where $u_d$ and $v_d$ are the 
least non-negative integers congruent to $u$ and $v$ 
respectively, modulo $d$.  So $\ce_d$ is  locally free at $p$ if 
and only if $d$ divides $u$ (and hence $v$).  If $u$ and $v$ are 
relatively prime, then no $\ce_d$ is locally free, and thus all 
$\ce_d$ are completely determined (up to isomorphism) by $(\ce_r, 
c_{r,1})$ (or $\pi^{\natural}\ce_r$) by 
$$\pi^{\natural}\ce_d := \pi^{\natural} \ce^{\otimes r/d}_r \otimes \co_{\tilde{X}}
(1/d(u-u_d)p^+ +1/d(v-v_d)p^-)$$ and 
$\ce_d=\pi_*\pi^{\natural}\ce_d$.  However, if $u$ and $v$ are 
not 
 relatively prime, but rather have $\gcd(u,v)=\ell>1$, then 
the root $\ce_\ell$ is locally free at $p$, and hence requires 
the additional gluing datum of an $\ell$th root of  unity 
(non-canonically determined) to construct $\ce_\ell$ from 
$\pi^{\natural}\ce_\ell$.  The remainder of the spin structure 
can clearly be reconstructed from the two pieces 
$c_{r,\ell}:\ce^{\otimes r/\ell}_r \rTo \ce_\ell$ and 
$c_{\ell,1}:\ce^{\otimes \ell}_\ell \rTo \omega.$ 

When the spin structure has no singularity (i.e., $\ce_r$ is 
locally free) at a node of the underlying curve, this corresponds 
to E.~Witten's definition of the {\em Ramond sector} of 
topological gravity 
\cite{witten:r-spin-conj,witten:N-matrix-model}; whereas, when 
$u$ and $v$ are non-zero, the spin structure is what Witten calls 
a {\em generalized Neveu-Schwarz sector}. If $\gcd(u,v)=\ell>1$ 
then we sometimes say that the spin structure is {\em 
semi-Ramond}, corresponding to the fact that $\ce_{\ell}$ is 
locally free (Ramond) but $\ce_r$ is not.

\subsubsection{Examples}\label{222}

It is useful to consider a few examples of stable spin curves.  
Both of the examples in this section are relevant to the study of 
the  Picard group of the stack and will be important in 
Section~\ref{div-basic}. 

{\bf Example 1: Two irreducible components and one node.}\label{two-comp}

First consider a stable curve $X$ with two smooth, irreducible
components $C$ and $D$, of genus $k$ and $g-k$ respectively, 
meeting in one double point $p$.  In this case there exists a 
unique choice of $u$ and $v$ that makes the degree of $\pi^* 
\omega_X (-u p^{+} - v p^{-}) \cong \omega_{\tilde X}(-(u-1)p^{+} 
- (v-1)p^{-}) $ divisible by $r$ on both components.  The 
resulting $r$th roots are locally free (Ramond) if and only if 
$\omega_X$ has an $r$th root, which is to say, if and only if $u$ 
and $v$ can be chosen to be $0$, or $2k-1 \equiv 0 
\pmod r$.  If $2k-1 \not\equiv 0 \pmod r$ then the resulting 
(Neveu-Schwarz) $r$th root corresponds to an $r$th root of 
$\omega_C(-(u-1)p^{+})$ on $C$ and an $r$th root of 
$\omega_D(-(v-1) p^{-})$ on $D$.  If $\gcd(2k-1,r)=1$ then all of 
the $d$th roots are Neveu-Schwarz.  Even if $\gcd(2k-1,r)=\ell$ 
is bigger than $1$, since the dual graph of $X$ is a tree, all 
gluing data for constructing $\ce_\ell$ from 
$\pi^{\natural}\ce_\ell$ will yield (non-canonically again) 
isomorphic $\ce_\ell$'s and hence isomorphic  nets of roots, and  
the spin curve corresponds to an 
 element in ${\stack S}_{k, 1}^{1/r,
  u-1} \cross {\stack S}_{g-k, 1}^{1/r, v-1}$.  Thus spin curves obey something like the splitting axiom of 
quantum cohomology (see \cite[\S4.1]{jkv} for more details). 

{\bf Example 2: One irreducible component and one
node.}\label{one-comp}

The second example is given by an irreducible stable curve $X$ with
one node.  In this case there are $r$ different choices of $u$ and $v$
that permit spin structures: either $u=v=0$, in which case the
resulting spin structure is locally free (Ramond); or $u \in
\{1,\dots, r-1\}$ and $v=r-u$, in which case the $r$th root
structure is not locally free (Neveu-Schwarz).  In this second 
case, if $\gcd(u,r)=1$, then the $r$-spin structure on $X$ is 
induced by an $r$th root of the bundle $\omega_{\tilde 
X}(-(u-1)p^{+} - (v-1)p^{-})$, and thus corresponds to an element 
of ${\stack S}_{g-1, 2}^{1/r, (u-1,v-1)}$. Since the points $p^+$ 
and $p^-$ of the normalization are not ordered, this gives a 
degree-$2$ morphism $\overline{\stack{S}}_{g,r}^{1/r,(n-1,v-1)} 
\rTo 
\sgrbar$.  But if $\ell=\gcd(u,r)>1$, then for a given $\pi^{\natural}\ce_\ell$ 
there are $\ell$ distinct choices of gluing data, and $\ell$ 
choices of  $\ce_\ell$.  This gives $\ell/2$ distinct morphisms 
$\overline{\stack{S}}_{g,2}^{1/r,(u-1),(v-1)}$ into $\sgrbar$.

\subsection{Properties of the Stack of Spin 
Curves}\label{moduli}\label{geometry}
\subsubsection{Basic Properties}

The main facts we need to know from \cite{jarvis:geometry} about 
the stack of stable spin curves are contained in the following 
theorem.  
 
\begin{theorem}\cite[Thms 2.4.4, and 
3.3.1]{jarvis:geometry})\label{finite}
 The stack 
$\overline{\stack{S}}^{1/r,\mathbf{m}}_{g,n}$ of $n$-pointed 
$r$-spin curves of genus $g$ and type $\mathbf{m}$ is a smooth, 
proper, Deligne-Mumford stack over $\mathbb{Z}[1/r]$, and the 
natural forgetful morphism $\overline{\stack{S}}^{1/r,\bm}_{g,n} 
\rTo 
\stackmgnbar$ is  finite and surjective. 

Moreover, if $\ell_{g,r}(\bm)$ is defined as 
$$\ell_{g,r}(\bm):=\left\{ \begin{array}
{ll} 1 &\text{ if } g=0 \text{ and } r|2+\sum m_i\\ 
\gcd(r,m_1, \dots, m_n) &\text{ if } g=1 \text{ and } r|\sum m_i\\
\gcd (2,r,m_1,\dots, m_n) &\text{ if } g>1 \text{ and } r| \sum m_i +2-2g\\
0 &\text{ otherwise}\\  
\end{array} \right.$$
and if $d_{g,r}(\bm)$ is defined to be the number of (positive) 
divisors of $\ell_{g,r}(\bm)$ (including $1$ and 
$\ell_{g,r}(\bm)$ itself), then 
$\overline{\stack{S}}^{1/r,\bm}_{g,n}$ is the disjoint union of 
$d_{g,r}(\bm)$ irreducible components. 

Also, $\overline{\stack{S}}^{1/r,\bm}_{g,n}$ contains the stack 
of smooth spin curves $\stacksgrnm$ as an open dense substack, 
and the coarse moduli spaces of $\sgrnmbar$ and $\sgrnm$ are 
normal and projective (respectively, quasi-projective). 
\end{theorem}

Although the actual details of the definition of a family of 
$r$-spin curves are not necessary for this paper, we do need the 
description from \cite{jarvis:geometry} of the universal 
deformation of a stable spin curve.  

\begin{theorem}(\cite[Thm. 2.4.2]{jarvis:geometry}):\label{univdef} 
Given a stable spin curve $(X, \{\ce_d,c_{d,d'}\})$ with 
singularities $q_i$ of order $\{u_i,v_i\}$, the universal 
deformation space of $(X, \{\ce_d, c_{d,d'}\})$ is the cover $$ 
\spec \mathfrak{o}[[\tau_1, \dots, \tau_m, t_{m+1}, \dots, t_{3g-3+n}]] \rTo \spec \mathfrak{o} [[t_1, \dots, t_{3g-3+n}]],$$
where $t_i=\tau^{r_i}_i$, and $r_i=r/\gcd (u_i,v_i)$, the scheme 
$\spec \mathfrak{o}[[t_1, \dots, t_{3g-3+n}]]$ is the universal 
deformation space for the underlying curve $X$, and the nodes 
$q_1, \dots, q_m$ correspond to the loci of vanishing of $t_1, 
\dots, t_m$, respectively. 
\end{theorem}

 \subsubsection{Relations Between the Different Stacks}\label{13}

There are several natural morphisms between the stacks.
\begin{enumerate}
\item There is a canonical isomorphism from $\stacksgrnmbar$ to 
$\overline{\stack{S}}_{g,n}^{1/r,\bm'}$ where $\bm'$ is an 
$n$-tuple whose entries are all congruent to $\bm \mod r$; namely 
for any net $\{\ce_d, c_{d,d'}\}$ of type $\bm$, let $\{\ce'_d, 
c'_{d,d'}\}$ be the net given by $\ce'_d =\ce_d \otimes \co(1/d 
\sum (m_i-m'_i)p_i)$ where $p_i$ is the $i$th marked point, and 
$c'_{d,d'}$ is the obvious homomorphism.  Because of this 
canonical isomorphism, we will often assume that all the $m_i$ 
lie between $0$ and $r-1$ (inclusive). 
\item There's a morphism $\overline{\stack{S}}_{g,n+1}^{1/r,\bm'}
\rTo^{\pi} \sgrnmbar$, where $\bm'$ is the $(n+1)$-tuple $(m_1, \dots, m_n, 0)$; and $\pi$ is the morphism which
simply forgets the $(n+1)$st marked point.  If $m_{n+1}$ is not 
congruent to zero $\mod r$, the degree of $\omega(\bm')$ is 
$2g-2-\sum^{n+1}_{i=1}m_i$ and the degree of $\omega(\bm)$ is 
$2g-2-\sum^n_1 m_i$.  Since both cannot be simultaneously 
divisible by $r$ at least one stack is empty, and there is no 
morphism.  This morphism, sometimes called ``forgetting tails,'' 
is not the universal curve over 
$\overline{\stack{S}}^{1/r,\bm}_{g,n}$, although it is birational  
to the universal curve.  In particular, the two are isomorphic 
over the open stack $\stack{S}^{1/r,\bm}_{g,n}$. 
\item 
If $s$ divides $r$, then there is a natural map $ [r/s]: 
\sgrnmbar 
\rTo \overline{\stack{S}}^{1/s,\bm}_{g,n}$, which forgets all of the roots and homomorphisms 
in the $r$th-root net except those associated to divisors of $s$.
\end{enumerate}

\section{Picard group}\label{2}
        \subsection{Definitions}\label{21}
     
Throughout this section we will assume that $g\geq 2$ and $n=0$, 
or that $g=1=n$.  
        
By the term {\em Picard group} we mean the Picard group of the
moduli functor; that is to say, the Picard group is the group of 
line bundles on the stack.  By a {\em line bundle $\cl$ on the 
stack $\stacksgrbar$}, we mean a functor that takes any         
family of spin curves $\mathfrak{Q} = (\cx/S, 
\{\ce_d,c_{d,d'}\})$ in $\stacksgrbar$ and assigns to it a line 
bundle $\cl(\mathfrak Q)$ on the scheme $S$, and which takes any 
morphism of spin curves $f:\mathfrak Q / S \rightarrow 
\mathfrak P /T$ and assigns to 
        it an isomorphism of line bundles $\cl(f): \cl(\mathfrak Q)
        \irightarrow f^{*} \cl(\mathfrak P)$, with the condition that
        the isomorphism must satisfy the cocycle condition (i.e., the
        isomorphism induced by a composition of maps agrees with the
        composition of the induced isomorphisms).  The groups $\pic
        \stacksgr$, $\pic \stackmg$, and $\pic \stackmgbar$ are defined similarly.
        For more details on Picard groups of moduli problems see
        \cite[pg.~50]{harris-mumford} or \cite[\S 5]{mumford:picard}.

 \subsection{Basic Divisors and Relations}\label{div-basic}
  \subsubsection{The Tautological Bundles}\label{221}
 
  Recall that for any family of stable curves $\pi
  :\cx \rightarrow S$, the Hodge class $\lambda (\cx/S)$ in $\pic \stack{M}_g$ is the determinant
  of the Hodge bundle (the push-forward of the canonical bundle)
  $$\lambda(\cx/S) := \det \pi_! \omega_{\cx/S} = \wedge^g \pi_{*}
  \omega_{\cx/S}.$$ It is well-known that $\pic \stackmg$ is the free Abelian
  group generated by $\lambda$ when $g>1$ (c.f. \cite{arbarello-cornalba:mg-pic}) and 
  for $g=1$ and $g=2$ $\pic \stack{M}_g$ is also generated by $\lambda$, but is cyclic of order $12$,
and $10$, respectively \cite[\S5.4] 
  {edidin-graham}.
  
  In a similar way we define a bundle $\mu=\mu^{1/r}$ in $\pic \stacksgrbar$ as the
  determinant of the $r$th root bundle.  In particular, if
  $\mathfrak Q = (\pi :\cx\rightarrow S,\{\ce_d,c_{d,d'}\})$ is a stable spin curve,
  then $\mu(\mathfrak Q) := \det \pi_! \ce_r = (\det R^0 \pi_* \ce_r)
  \tensor (\det R^1 \pi_* \ce_r)^{-1}$ on $S$.  Similarly, define 
  $\mu^{1/d}:= \det \pi_!\ce_d$ for each $d$ dividing $r$.  Note 
  that the pullback of $\mu$ from $\stacksgsbar$ to 
  $\stacksgrbar$ via $[r/s]^*$ is exactly $\mu^{1/s}$.

 \subsubsection{Boundary Divisors Induced from $\stackmgbar$}\label{bdry-div}

In addition to $\lambda$ and $\mu$, there are elements of $\pic
\stacksgrbar$ that arise from the boundary divisors of $\stacksgrbar$.  Recall
that the boundary of $\stackmgbar$ consists of the divisors 
$\delta_i$ where $i\in \{0,\dots,\lfloor g/2 \rfloor\}$.  Here, 
when $i$ is greater than zero, $\delta_i$ is the closure of the 
locus of points in $\stackmgbar$ corresponding to stable curves 
with exactly one node and two irreducible components, one of 
genus $i$ and the other of genus $g-i$. And when $i$ is zero, 
$\delta_{0}$ is the closure of the locus of points corresponding 
to irreducible curves with a single node. 

As we saw in  Example~1 of Section~\ref{two-comp}, for any curve 
$X$ with exactly one node and two irreducible components of 
genera $i$ and $g-i$, there is a unique choice of integer $u(i)$ 
between $0$ and $r-1$ that determines a bundle $\cw$ whose $r$th 
roots define the spin structures on $X$.  If $2i\equiv 1 
\pmod r$, then $u(i)=v(i)=0$ and the bundle $\cw$ is just the 
canonical bundle  $\cw = 
\omega_X$.  In this case all the spin structures on $X$ are locally
free.  If on the other hand, $2i\not\equiv 1 \pmod r$, then there is a
unique choice of $u(i)$ with $r>u(i)>0$, and such that $2i-1-u(i)
\equiv 0 \pmod r$.  In this case  $v(i)=r-u(i)$, and $\cw$ is
not a line bundle on $X$, but rather a line bundle on the
normalization $\nu: \tilde X \rightarrow X$ at the node $q$.  If
$\nu^{-1}(q) = \{ q^{+},q^{-} \}$, then 
$$\cw = \omega_{\tilde{X}}((1-u(i)) q^{+}+(1- v(i)) q^{-}).$$ In this case the
spin structures on $X$ are constructed by pushing the $r$th roots 
of $\cw$ on $\tilde{X}$ down to $X$, i.e., using the $r$th roots 
of $\cw$ on the irreducible components to make a torsion-free 
sheaf on the curve $X$.  By Theorem~\ref{finite}, if $i$ and 
$g-i$ are at least $2$, and $r$ is odd, there is one irreducible 
divisor $\alpha_i$ lying above $\delta_i$.  And when $r$ is even, 
there are four divisors lying over $\delta_i$ (even-even, 
even-odd, odd-even, odd-odd). 

In general, let $D_{k,r}(\mathbf{m})$ denote the set of positive 
divisors of $\ell_{k,r}(\mathbf{m})$, as in Theorem \ref{finite}.  
These divisors index the set of irreducible components of 
$\overline{\stack{S}}^{1/r,\bm}_{k,n}$.  For each $i \geq 1$, and 
each $a 
\in D_{i,r}(u(i)-1)$ and $b \in D_{g-i,r}(v(i)-1)$, let 
$\alpha^{(a,b)}_i$ denote the irreducible divisor consisting of 
the locus of spin curves lying over $\delta_i$ with a spin 
structure of index $a$ on the genus-$i$ component and of index 
$b$ on the genus-$(g-i)$ component. 

Over a stable curve in $\delta_{0}$ there are also several 
choices of spin structure.  Indeed, as we saw in Example~2 of 
Section~\ref{one-comp}, for any choice of order $\{u,v\}$, there 
is again a unique bundle $\cw_u = \omega_{\tilde X}((1-u) 
q^{+}+(1-v) q^{-})$ so that any $r$th root of $\cw_u$ gives the 
pullback $\pi^{\natural}\ce_r$ to $\tilde{X}$, and thus the 
entire spin structure \emph{except for glue}.  And conversely, 
any spin structure comes from an $r$th root of some $\cw_u$ for 
$0 
\le u < r$ and a choice of glue (where the glue corresponds to an 
isomorphism $\eta:\ce_{\ell}|_{q^+} \irightarrow 
\ce_{\ell}|_{q^-}$ and $\ell=\gcd(u,v,r)$). Again 
for a particular choice of order $\{u,v\}$ of index $a \in 
D_{g-1,r}(u-1,v-1)$, and a particular choice of glue $\eta$, the 
corresponding divisor of spin curves of the given order, index 
and glue is irreducible. We denote these divisors by 
$\gamma_{j,\eta}^{(a)}$, where $j$ is the smaller of $u$ and $v$, 
$\eta$ is the gluing datum, and $a$ is the index.  Of course, 
since the two points of the normalization are not 
distinguishable, we have 
$\gamma^{(a)}_{j,\eta}=\gamma^{(a)}_{r-j,\eta^{-1}}$.  We denote 
the set of gluings for a given choice of $j$ by $\mathfrak{g}_j$,  
and we denote the set of gluings for a fixed $j$, modulo the 
equivalence $(r/2,\eta) \sim (r/2,\eta^{-1})$, by 
$\mathfrak{g}_j/S_2$.  Let $\alpha_i$ denote the sum 
$$\alpha_i:=\sum_{D_i(u(i)-1)\cross D_{g-i}(v(i)-1)} 
\alpha_i^{(a,b)}$$ and let $\gamma_j$ denote the sum $$\gamma_j:=\sum_{\substack{\eta\in \mathfrak{g}_j/S_2\\a \in D_{g-1} (j-1,r-j-1)}} \gamma^{(a)}_{j,\eta}.$$ 

Pullback along the forgetful map $p:\stacksgrbar 
\rightarrow \stackmgbar$ induces a map $p^*$ from $\pic \stackmgbar$ to 
$\pic \stacksgrbar$.  And it is clear that the image of the 
boundary divisors (which we will denote by $\delta_i$, regardless 
of whether we are working in $\pic 
\stackmgbar$ or $\pic \stacksgrbar$) is still supported on the boundary
of $\stacksgrbar$.  Indeed, for $i$ greater than one, $\delta_i$ 
is a linear combination of the $\alpha_i^{(a,b)}$'s, and  
$\delta_{0}$ is a linear combination of the 
$\gamma^{(a)}_{j,\eta}$'s. 

\begin{proposition}\label{boundary-prop}
Let $c_i:= \gcd(2i-1,r)= \gcd(u(i),v(i))$.  If $i \ge 1 $  then 
$$\delta_i = (r/c_i)\alpha_i $$ 
And if $d_j := \gcd(j,r)$, then we have
 $$\delta_{0} =
                         \sum_{0 \le j \le r/2} (r/d_j) \gamma_j.$$
\end{proposition}

\begin{proof}
  
  The universal
  deformation of a spin curve with underlying stable curve in
  $\delta_i$ is dependent only upon $u$ and $v$; in particular, if
  $c:=\gcd(u,v)$ and $r' := r/c$ then the forgetful map from the
  universal deformation of the spin curve to the universal deformation
  of the underlying curve is of the form $\spec \mathfrak o [[s]]
  \rightarrow \spec \mathfrak o [[s^{r'}]] $ (see Theorem \ref{univdef}).  Moreover, for $i\ge 1$ the morphism 
  $\stacksgrbar
  \rightarrow \stackmgbar$ is ramified along $\alpha_i$  to
  order $r' = r/c_i$.  Similarly, over $\delta_{0}$, the map $\sgrbar \rightarrow \mgbar$ is
  ramified along $\gamma_j$ to order $r' = r/d_j$.
  In either case the proposition follows.

\end{proof}

\subsubsection{Boundary Divisors Induced from Other-order Spin 
Curves} If $s$ divides $r$, say $sd=r$, then the natural map 
$[d]:\stacksgrbar \rTo \overline{\stack{S}}_{g}^{1/s}$ induces a 
homomorphism $[d]^*:\pic \overline{\stack{S}}_{g}^{1/s} \rTo \pic 
\stacksgrbar$.  Let $\alpha^{1/s}_{i}$, 
and $\gamma^{1/s}_{j}$ also indicate the images of the  
corresponding boundary divisors under the map $[d]^*$.

\begin{proposition}
  Let $u(i)$ indicate the unique integer $0 \le u(i) < r$ such that
  $2i-1-u(i) \equiv 0 \pmod r$. For each divisor $s$ of $r$, let $c^{1/s}_{i}$ be $\gcd
  (u(i),v(i),s)=\gcd(u(i),s)$ and let $d^{1/s}_{i}$ be $\gcd(j,(r-j),s)= \gcd(j,s)$.  In $\pic \stacksgrbar$ the boundary
  divisors $\alpha^{1/s}_{i}$, and $\gamma^{1/s}_{i}$ are given
  in terms of the elements $\alpha^{1/r}_{i}$, and
  $\gamma^{1/r}_{j}$ as follows.

  $$\alpha^{1/s}_{i} = \frac{r}{c^{1/r}_{i}} \frac{c^{1/s}_{i}}{s} 
  \alpha^{1/r}_{i}$$
  
  And $\gamma^{1/s}_{k}$ behaves similarly, but now the index $j$ behaves like
  $u \pmod s$, and there are several choices of $k$ between $0$ and
  $r/2$ 
  such that $k \equiv \pm j \pmod s$ ($\gamma^{1/r}_{k}$ is the same divisor as
  $\gamma^{1/r}_{r-k}$).  This gives several divisors $\gamma^{1/r}_{k}$ over each $\gamma^{1/s}_{j}$:
 
$$\gamma^{1/s}_{j} = \sum_{\substack{k \equiv \pm j (\operatorname{mod} s)\\0 \leq k \leq r/2}} \frac{r}{d^{1/r}_{k}} 
\frac{d^{1/s}_{k}}{s} \gamma^{1/r}_{k}$$

\end{proposition}
                    
\begin{proof}
  The proof is straightforward and very similar to the proof of Proposition~\ref{boundary-prop}.  Only note that
  $\gamma^{1/r}_{k}$ lies over $\gamma^{1/s}_{j}$ if  $k \equiv j
  \pmod s$ {\em or if}  $k \equiv - j \pmod s$.
\end{proof}

\subsubsection{Independence of Some Elements of $\pic \stacksgrbar$}\label{224}

 The following two results are generalizations of Cornalba's results
 \cite[Proposition 7.2]{cornalba:theta}.  Their
 proofs are essentially the same as those in \cite{cornalba:theta},
 and so will just be sketched here.

\begin{proposition}
The forgetful map  $\stacksgrbar
\rightarrow \stackmgbar$ induces an injective homomorphism  on the Picard groups.
\end{proposition}

\begin{proof}
  The method of proof is simply to recall
  from \cite{arbarello-cornalba:mg-pic} that there are families of
  stable curves $\cx/S$, with $S$ a smooth and complete curve, such
  that the vectors $(\deg_S \lambda,\deg_S \delta_{0}, \dots, \deg_S
  \delta_{\lfloor g/2 \rfloor})$ are independent.  After suitable base
  change, one can install a spin structure on the families in
  question.  (This follows, for example, from the fact that, 
  after base change, a spin structure can be installed on the 
  generic fibre of $\cx/S$, and since the stack is proper over 
  $\stackmgbar$, this structure can be extended--again after 
  possible base change--to the entire family).
   Since the effect of base change is to multiply the
  vectors' entries by a constant, the vectors are still independent,
  and thus the elements $\lambda$ and $\delta_i$ are all independent
  in $\pic \stacksgrbar$.

\end{proof}

J.~Koll\'ar pointed out to me the following alternate proof:
Note that the stacks in question are both smooth.  Thus the result
follows from the fact that for any finite cover $f:X' \rightarrow X$
of a normal variety $X$ and for any line bundle $\cl$ on $X$, the line
bundle $f_*f^* \cl$ is a multiple of $\cl$ (c.f.,
\cite[6.5.3.2]{ega2}). In other words, for some $n$ we have $$f_*f^*
\cl = \cl^{\tensor n}.$$
But since the Picard group of $\stackmgbar$ has no
torsion, the pullback of any line bundle to $\stacksgrbar$ cannot be
trivial.

\begin{proposition}\label{lag-indep}
The elements $\lambda$, $\alpha^{(a,b)}_i$ and 
$\gamma^{(a)}_{j,\rho}$ are independent in $\pic 
\stacksgrbar$ when $g>1$. 
\end{proposition}

\begin{proof}

Again the idea is to install a spin structure on a family of 
curves $\cx/S$ over a smooth, complete curve, where the degrees 
of $\lambda$, and the $\delta_i$ are known.  In particular, given 
two curves $S$ and $T$ of genera $i$ and $g-i$, respectively, fix 
$t \in T$ and let $s$ vary in $S$.  Consider the family $\cx / S$ 
constructed by joining the two curves at the points $s$ and $t$ 
\cite[\S 7]{harris-mumford}. Then the degrees of $\lambda$ and 
$\delta_j$ are all zero on $S$, except when $j = i$, and then 
$\deg_S \delta_i = 2-2i$.  We can construct an $r$th root of 
$\omega_{S}((1-u(i))s)$ on $S$ and an $r$th root of $\omega_T 
((1-v(i))t)$ on $T$, and thus an $r$-spin structure on $\cx/S$.  
Moreover, the two $r$th roots can be chosen to be of any index in 
$D_{i,r}(u(i)-1)$ or $D_{g-i,r}(v(i)-1)$, respectively; and thus 
$\cx/S$ can be endowed with a spin structure of any index $(a,b)$ 
in $D_{i,r}(u(i)-i)\cross$ $D_{g-i,r}(v(i)-1)$ along every fibre.  
Because the $\alpha^{(a,b)}_i$ are disjoint, the degree of 
$\alpha^{(a',b')}_i$ for any other index $(a',b')$ must be zero 
on $S$, but $\deg_S 
\delta_i \neq 0$ implies that $ \deg(\sum_{(a,b)}\alpha^{(a,b)}_i)$ is non-zero.

Consequently, in any relation of the form $0 = \ell \lambda + 
\sum e^{(a,b)}_i \alpha^{(a,b)}_i + \sum c^{(a)}_{j,\rho} \gamma^{(a)}_{j,\rho}$, the coefficients $e_i$ must 
all vanish.  And thus a relation must be of the form $ 0 = \ell 
\lambda + \sum c^{(a)}_{j,\rho} \gamma^{(a)}_{j,\rho}$.  But a similar method shows that 
the coefficients $c^{(a)}_{j,\rho}$ must also be zero.  In 
particular, consider the family $\cy/C$ constructed by taking a 
general curve $C$ of genus $g-1$ and identifying one fixed point 
$p$ with another, variable point $q$ \cite[\S 7]{harris-mumford}.  
Again one may produce an $r$-spin structure on $\cy/C$ of any  
type. And $\deg_C \lambda = 
\deg_C \gamma_{k,\psi}^{(b)}  = 0 $ whenever $\gamma^{(b)}_{k,\psi}\neq \gamma^{(a)}_{j,\rho}
$, but $\deg_C 
\delta_{0} = \sum_{j,\rho,a} r/d(j) \deg \gamma^{(a)}_{j,\rho}$ is equal to $2-2g$, which is non-zero (since $g>1$).  Thus $c^{(a)}_{j,\rho} =
0$, and so also $\ell = 0$.
\end{proof}

\subsection{Less-Obvious Relations and Their Consequences}\label{23}
  
Another important question about these bundles is what relations exist
between the bundles $\mu,$ $\lambda,$ and the various boundary
divisors.  In this section we provide a partial answer.

\subsubsection{Main Relation}  The main result of this section 
is Theorem~\ref{mainthm}, which provides relations between 
$\mu^{1/r}$, $\mu^{1/s}$, $\lambda$, and some boundary divisors.  
The proof will be given in Section ~\ref{aux-rel}.  To state the 
result we need the following definitions.  Motivation for the 
somewhat-peculiar notation will be made clear in the following 
section. 

\begin{defn}\label{unique}
  Let $u(i)$ denote, as in the previous section, the unique integer
  $0\le u(i) < r$ such that $2i-1 \equiv u(i) \pmod r$. Let $v(i) =
  r-u(i)$, $c^{1/s}_{i}:=\gcd(u(i),v(i),s) = \gcd(2i-1,s)$, and $d^{1/s}_{j}:=\gcd(j,r-j,s)= \gcd(j,s)$ and 
  let $c_i = c^{1/r}_{i},d_j=d^{1/r}_{j}$.
  Then we define $<\tilde \ce_r, \mathfrak E_{r}>$ to be the element in $\pic
  \stacksgrbar$ defined by the following boundary divisor: $$<\tilde
  \ce_r, \mathfrak E_{r}> := \sum_{1 \le i \le g/2} (u(i) v(i)/c_i) \alpha_{i} 
   + \sum_{0\le j \le r/2} (j (r-j)/d_j) \gamma_{j}$$
\end{defn}

Similarly, for $s$ dividing $r$ we make the following 
definitions. 

\begin{defn}
 Let $u'(i)$ denote the unique integer $0\le u'(i) < s$
such that $2i-1 \equiv u'(i) \pmod s$. Let $v'(i) = s-u'(i)$.  
And let $<\tilde \ce^{1/s}, \mathfrak E_{s}>$ denote the  
following boundary divisor in $\pic 
\stacksgrbar$: 

\begin{eqnarray*}
\displaystyle 
<\tilde \ce^{1/s}, \mathfrak E_{s}>& = &\sum_{1 
\le i \le g/2} 
\frac{u'(i) v'(i)}{c^{1/s}_{i}} \alpha^{1/s}_{i}  
 + \sum_{0\le j \le s/2} \frac{j (s-j)}{d^{1/s}_{j}} \gamma^{1/s}_{j}\\  
&=&\frac{r}{s} \sum_{1 \le i \le g/2} \frac{(u'(i) 
v'(i))}{c^{1/r}_{i}} \alpha^{1/r}_{i} + \frac{r}{s} 
\sum_{\substack{1\le j \le s/2\\ 1\le k \le r/2 \\ k \equiv \pm j
\pmod s}} \frac{(j (s-j))}{d^{1/r}_{k}} \gamma^{1/r}_{k} 
\end{eqnarray*} 
\end{defn}

\begin{defn}
 Let $\delta$ indicate the element of $\pic \stacksgrbar$ defined by the
sum of the boundary divisors pulled back from $\pic \stackmgbar$:
\begin{align*} \delta & = \sum_{i=0}^{g/2} \delta_i\\ & = 
 \sum_{0 \le j\le r/2} (r / d_j) \gamma_{j} + \sum_{1\le i \le g/2}
(r / c_i) \alpha_{i} \\ & =\sum_{0 \leq k \leq s/2} 
\left(\frac{s}{d^{1/s}_{k}} \right) \gamma^{1/s}_{k} + \sum_{1 \leq i 
\leq g/2} \left( \frac{s}{c^{1/s}_{i}}\right)\alpha^{1/s}_{i} \\
\end{align*} 
\end{defn}

With all of this notation in place we can now state the main theorem.

\begin{theorem}\label{mainthm}
  In terms of the notation defined above, the following relations hold
  in $\pic \stacksgrbar$:
\begin{align*}
  r <\tilde \ce_r, \mathfrak E_{r}> &= (2r^2 - 12 r + 12)\lambda - 2r^2 \mu +
  (r-1)
  \delta  \\
  s <\tilde \ce_s, \mathfrak E_{s}> &= (2s^2 - 12 s + 12)\lambda - 2s^2 \mu_s +
  (s-1)
  \delta \\
\end{align*}
\end{theorem}

This may also be written in terms of $\alpha$ and $\gamma$ 
instead of  $<\tilde \ce_r, \mathfrak E_{r}>$, $<\tilde 
\ce_s, \mathfrak E_{s}>$, and $\delta$.  But the final 
expression is greatly simplified if we define 
$$\sigma^{1/s}_{k} 
:= \sum_{\substack{1\le i \le g/2 \\ i \equiv k \pmod s}} 
\alpha^{1/s}_{i}.$$ Let $\sigma_k =\sigma^{1/r}_{k}$, 
and if $k$ is not an integer, then let $\sigma^{1/s}_{k}=0$. 

Also, we will continue to use the notation  
$c^{1/s}_{j}=\gcd(2i-1,s),$  $c_i=c^{1/r}_{i}$, 
$d^{1/s}_{j}=\gcd(j,s)$, and $d_j=d^{1/r}_{j}$. 

\newtheorem*{thmbis}{\bf Theorem~\ref{mainthm}.bis}
\begin{thmbis}
In terms of the  notation defined above, the following relation
holds in $\pic \stacksgrbar$:
\begin{align*}
(2r^2 - 12 r + 12)\lambda - 2r^2 \mu   
& = (1-r) (\sigma_{\frac{r+1}{2}} + \gamma_0)\\
  &+  \sum_{1 < k < \frac{r}{2}} 2 (r/c_k) (r k - 2 k^2  + 2k - r) \sigma_k \\
  &+ \sum_{\frac{r}{2} + 1 < k < r}  2(r/c_k)
(3r k - 2k^2 + 2k  - 2r - r^2) \sigma_k\\
  &+ \sum_{1 < j \le r/2} (r/d_j) ( j (r-j) - (r - 1) ) \gamma_j   \\
\end{align*}

Similarly,
\begin{align*}
(2s^2 - 12 s + 12)\lambda - 2s^2 \mu^{1/s}   & = (1-s) 
(\sigma_{\frac{s+1}{2}} + \gamma^{1/s}_{0})\\ 
  &+  \sum_{1 < k < \frac{s}{2}} 2 (s/c^{1/s}_{k}) (s k - 2 k^2  + 2k - s) \sigma^{1/s}_{k} \\
  &+ \sum_{\frac{s}{2} + 1 < k < s}  2(s/c^{1/s}_{k})
(3s k - 2k^2 + 2k  - 2s - s^2) \sigma^{1/s}_{k}\\
  &+ \sum_{1 < j \le s/2} (s/d^{1/s}_{j}) ( j (s-j) - (s - 1) ) \gamma^{1/s}_{j}   \\
\end{align*}
\end{thmbis}

The proof of Theorem~\ref{mainthm} will be postponed until the next
section.  The proof of Theorem~\ref{mainthm}.bis is a
straightforward, but tedious calculation, accomplished simply by
writing out the definitions of all of the different terms, and
applying Theorem~\ref{mainthm}.

The main thing to note about the relations of Theorems~\ref{mainthm}
and \ref{mainthm}.bis is that they hold in $\pic \stacksgrbar$ and not just
in $\pic \stacksgrbar \tensor \mathbb Q$---that is to say, not just modulo
torsion.

We also have the following immediate corollaries:

\begin{corollary}
  If $s$ divides $r$  the following relations hold in $\pic \stacksgr$:
\begin{align*}
  (2 r^2 - 12 r + 12)\lambda &= 2r^2 \mu^{1/r}\\
  (2 s^2 - 12 s + 12)\lambda &= 2s^2 \mu^{1/s}\\
 2 r^2(s^2-6s+6)\mu^{1/r}  &=  2s^2(r^2-6r+6) \mu^{1/s}  
\end{align*}
\end{corollary}

\begin{corollary}\label{table}
 The following special cases of the relations in
Theorems~\ref{mainthm} and \ref{mainthm}.bis hold:

$$
\begin{array}{|l|rcl|}

\hline
r=2 & 4\lambda + 8 \mu^{1/2}   & = & \gamma_0 \\ 
\hline

r=3 & 6\lambda + 18 \mu^{1/3}  & = & 2 \gamma_0 + 2 \sigma_2 \\ 
\hline

r=4 & 4\lambda + 32 \mu^{1/4}  & = & 3 \gamma_{0} - 2\gamma_{2}\\ 
    & 4\lambda + 8 \mu^{1/2}   & = & \gamma_{0} + 2 \gamma_{2}\\ 
\hline

r=5 & 2\lambda - 50 \mu^{1/5}  & = & -4 \gamma_0 + 10 (\gamma_2 + 
\sigma_2 
                                + \sigma_4) - 4 \sigma_3\\
\hline
r=6 & 12\lambda - 72 \mu^{1/6} & = & -5 \gamma_{0} + 9 
\gamma_{2} + 
                                8 (\gamma_{3} + \sigma_{2} + \sigma_{5}) \\
    & 4\lambda + 8 \mu^{1/2}   & = & \gamma_{0} + 3\gamma_{2} \\
    & 6\lambda + 18\mu^{1/3}   & = & 2 \gamma_0 +4 \gamma_3 + 4(\sigma_2+\sigma_5) \\
\hline
r=7 & 26 \lambda -98 \mu^{1/7} & = & -6(\gamma_0+\sigma_4) + 28 
(\gamma_2 + \sigma_3 + \sigma_5) + 42(\gamma_3 + \sigma_2 + 
\sigma_6) \\
\hline
r=8 & 44 \lambda - 128 \mu^{1/8} & = & -7\gamma_0 + 20 \gamma_2 +  
18 \gamma_4 + 64(\gamma_3+\sigma_2+\sigma_3+\sigma_6+\sigma_7) 
\\ & 4 \lambda-32\mu^{1/4} & = & 3 \gamma_0 + 6 \gamma_4 - 4 
\gamma_2 
\\
& 4 \lambda -8\mu^{1/2} & = & \gamma_0 +4\gamma_2+2 \gamma_4\\
\hline

\end{array}
$$

\end{corollary}

Note that in the case when  $r=2$, the relation $4\lambda + 8
\mu^{1/2} = \gamma_0$  is exactly the content of Theorem~(3.6) in
\cite{cornalba:spin-pic}.  

\subsubsection{Proof of the Main Relation and Some Other Relations}
\label{aux-rel}

Although some of the proof of Theorem~\ref{mainthm} resembles that of
Theorem~(3.6) in \cite{cornalba:spin-pic}, the proof of
Theorem~\ref{mainthm} requires many involved calculations for the
case $r>2$ that do not arise when $r=2$.

To accomplish the proof without annihilating torsion elements, we
need some tools other than the usual Grothendieck-Riemann-Roch; in
particular, we'll use the following construction of Deligne
\cite{deligne:determinant}. 
\begin{defn} 
  Given two line bundles $\cl$ and $\cm$ on a family of semi-stable curves
  $f:\cx\rightarrow S$, define a new line bundle, the Deligne product
  $<\cl,\cm>$ on $S$, as
$$<\cl,\cm> := (\det f_! \cl\cm) \tensor (\det f_! \cl)^{-1} \tensor
(\det f_!  \cm)^{-1}\tensor \det f_! \co_\cx.$$
\end{defn}
Here, as earlier, $\det f_! \cl$ indicates the line bundle $(\det 
R^0 f_* \cl) \tensor (\det R^1f_* \cl)^{-1}$.  We will write this 
additively as $<\cl, \cm>=\det f_!\cl,\cm-\det f_! 
\cl-\det f_! \cm + \det f_! \co$.  Deligne shows that this operation is symmetric and 
bilinear.  It is also straightforward to check that when the base 
$S$ is a smooth curve then the degree on $S$ of the Deligne 
product $<\cl,\cm>$ is just the intersection number of $\cl$ with 
$\cm$; that is to say, if $(- . 
- )$ is the intersection pairing on $\cx$, then $$\deg_S 
<\cl,\cm> = (\cl . 
\cm).$$ 

Using this notation, and writing $\omega$ for the canonical bundle of $\cx/S$,
Serre duality takes the form
$$<\omega,\cm> = \det f_!(\cm ^{-1}) - (\det f_! \cm),$$
or 
$$\det f_! (\cl \omega) = \det f_!(\cl^{-1}).$$ 
Deligne proves the following variant of Grothendieck-Riemann-Roch 
that holds with integer, rather than just rational coefficients. 

\begin{lemma}(Deligne-Riemann-Roch)\label{GRR}
Let  $f:\cx\rightarrow S$ be a family of stable curves, and let $\cl$ be
any line bundle on $\cx/S$.  If $\omega$ denotes the canonical
sheaf of $\cx/S$, then 
$$ 2(\det f_! \cl) = <\cl,\cl> - <\cl,\omega>+2(\det
f_! \omega).$$
\end{lemma}

We also need some additional line bundles associated to a family 
of $r$-spin curves $(\cx/S,\{\ce_d,c_{d,d'}\})$ that intuitively 
amount to measuring the difference between $\ce_r$ and a ``real'' 
$r$th root of the canonical bundle.  To construct these bundles 
we first recall from \cite{jarvis:spin} some canonical 
constructions related to any $r$th root $(\ce_r,c_{r,1})$ on a 
curve $\cx/S$   

\begin{enumerate}
\item There is a semi-stable curve $\pi: \tilde \cx_{\ce_r} 
\rightarrow 
\cx$, uniquely determined by $\cx$ and $\ce_r$, so that the stable model of $\tilde \cx$ is 
$\cx$.
\item On $\tilde \cx$, 
there is a unique line bundle $\co_{\tilde \cx}(1)$, and a 
canonically-determined, injective map $\beta: \co_{\tilde 
\cx}(1)^{\tensor r} \rightarrow 
\omega_{\tilde \cx /S}$,
the push forward $\pi_* \co_{\tilde \cx}(1) $ is $\ce_r$, and the 
map $b$ is induced from $\beta$ by adjointness.  
\item Moreover, the 
degree of $\co_{\tilde \cx}(1)$ is one on any exceptional curve 
in any fibre. 
\end{enumerate}

To simplify, we will denote $\co_{\tilde \cx_{\ce_r}}(1)$ by 
$\tilde 
\ce_r$.  And we define a new bundle $\mathfrak E$ by
$$\mathfrak E_r = \omega_{\tilde \cx} \tensor \tilde \ce_r^{\tensor -r}.$$
The advantage of using $\mathfrak E_r$ is that it is completely
supported on the exceptional locus, and thus it is easy to 
describe explicitly; and its Deligne product with other line 
bundles on $\tilde 
\cx$ easy to compute.

Similarly, if $s$ divides $r$, then there is a family $( \ce_s, 
c_{s,1})$ of $s$th roots on $\cx$.  And corresponding to this 
family, there is a line bundle $\tilde 
\ce_s$ on $\tilde \cx_{\ce_s}$ which pushes down to $\ce_s$, and 
another line bundle $\mathfrak E_s = \omega_{\tilde \cx_s} 
\tensor (\tilde \ce ')^{\tensor -s}$. 

\begin{theorem} \label{335}\label{bleh} 

The Deligne products $<\tilde \ce_r,\mathfrak E_r>$ and 
$<\tilde{\ce_s},\mathfrak{E}_s>$ agree with Definition 
\ref{unique}.  Also, the Deligne product 
  of the canonical bundle $\omega$ with $\mathfrak E$ 
is trivial.  In particular, 
we have
\begin{enumerate}
\item The  Deligne product $<\tilde \ce_r, \mathfrak E_r>$ is exactly $$
  \sum_{1 \le i \le g/2} (u(i) v(i)/c_i) \alpha^{1/r}_{i}  +
  \sum_{0\le j \le r/2} (j (r-j)/d_j) \gamma^{1/r}_{j}.$$
\item The Deligne product $<\tilde \ce_s, \mathfrak E_s>$ is exactly
   $$ r/s \! \! \sum_{1 \le i \le g/2} (u''(i) v''(i)/c^{1/r}_i)
     \alpha^{1/r}_{i} + r/s \! \! \! \sum_{\substack{1\le j \le s/2\\ 1\le k \le
       r/2 \\ k \equiv j \pmod s}} j (s-j)/d^{1/r}_j \gamma_{k^{1/r}}. $$
\item \label{we} And $<\omega_{\tilde \cx_r}, \mathfrak E_r> = \co= <\omega_{\tilde{\cx}_s},\mathfrak E_s>$.

\end{enumerate}
\end{theorem}

\begin{proof}
  For $\cx$ smooth, $\mathfrak E_r$ 
  is canonically isomorphic to $\co$, and thus in the smooth case,
  $<\tilde \ce_r, \mathfrak E_r>$
  and  $<\omega_{\tilde \cx}, \mathfrak E_r>$
are trivial.

In the general case these products are all integral linear
combinations of boundary divisors.  To compute the coefficients 
we evaluate degrees on families of curves $\cx/S$ parameterized 
by a smooth curve $S$ and having smooth generic fibre.  In this 
case, since $\mathfrak E$ is supported on the exceptional locus 
of $\tilde \cx$, whereas the canonical bundle has degree zero on 
the exceptional locus, the product $<\omega_{\tilde 
\cx},\mathfrak E_r>$ must be trivial. 

Computing the coefficients of $<\tilde \ce_r, \mathfrak E_r>$ 
requires that we consider the local structure of $\mathfrak E_r$ 
near an exceptional curve.  Recall from \cite{jarvis:geometry} 
that for any singularity where $\ce$ has order $\{u,v\}$ with 
$u,v>0$, and such that if $c:=\gcd(u,v)$, $u':=u/c$, $v':=v/c$, 
and $r':=r/c$, the underlying singularity of $\cx/S$ is 
analytically isomorphic to $\spec (\hat{\co}_{S,s} [[x,y]]/{x 
  y-\tau^{r'}})$, where $\tau$ is an element of $\hat{\co}_{S,s}$. $\ce_r$ is generated by 
two elements, say $\nu$ and $\xi$, with the relations $x \nu = 
\tau^{v'} \xi$ and $y\xi = \tau^{u'} \nu$.  Moreover, over such a 
singularity, $\pi:\tilde \cx \rightarrow \cx$ is locally given as 
$$\pi:\proj_{A}(A[\nu ,\xi]/(\nu x-\tau^{v'} \xi, \nu \tau^{u'}-\xi y))
\rightarrow \spec A,$$ where
$A$ is the local ring of $\cx$ at the singularity.  The exceptional
curve, call it $D$, is defined by the vanishing of $x$ and $y$, and we
have a situation like that depicted in Figure~\ref{xtilde}.

\begin{figure}
\includegraphics{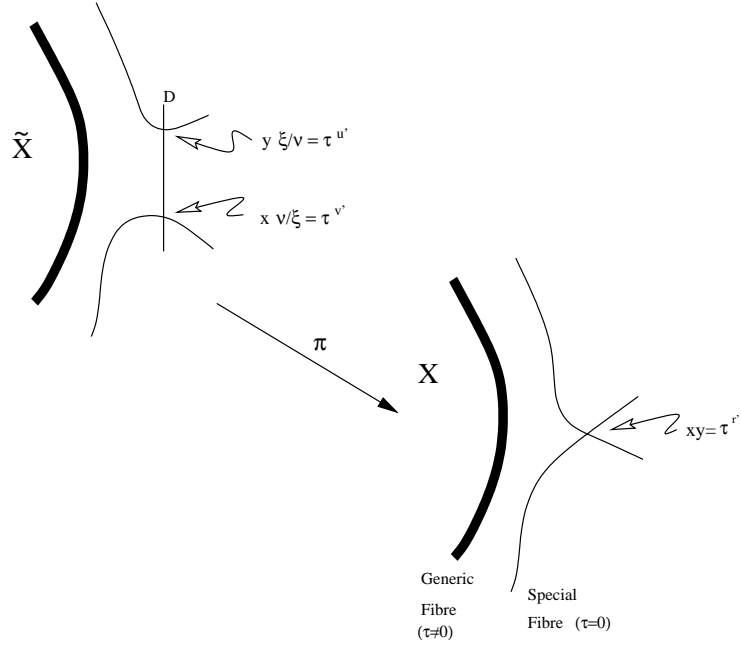}
\caption{\label{xtilde}
Local structure of $\pi: \tilde \cx \rightarrow \cx$ over a
singularity of $\ce$.}
\end{figure}

We need to express $\mathfrak E_r$ in terms of a divisor, but 
this is easy since it is supported completely on the exceptional 
locus. $\mathfrak E_r$ is locally of the form $\co_{\tilde \cx} 
(nD)$, for some integer $n$.  And any Weil divisor of the form 
$nD$ is Cartier if and only if $u'$ and $v'$ both divide $n$.  
Moreover, it is easy to see that if $nD$ is Cartier, then when 
restricted to the exceptional curve $D$, the degree of $nD$ is $ 
-n/u' - n/v'$.  Finally, since $\mathfrak E_r$ has degree $-r$ on 
$D$, we have $n = u' v' c = u v /c$ so that 
        $$\mathfrak E_r = \co_{\tilde \cx} ( (u v/c) D).$$

Now, to compute the coefficients of $<\tilde \ce_r, \mathfrak 
E_r>$ 
 just note that for families $f: \cx \rightarrow S$ over a smooth base
 curve $S$ with smooth generic fibre, the degree of $<\tilde \ce_r,
 \mathfrak E_r>$ is just the intersection number $(\tilde \ce_r,  
 \mathfrak E_r)$.  Thus if $D_p$ indicates the exceptional curve over a
 singularity $p$, and if $\{u_p,v_p\}$ indicates the order of $(\ce_r,c_{r,1})$ near $p$, then
\begin{eqnarray*}\deg_S <\tilde \ce_r, \mathfrak E_r> & = &
 \sum_{\substack{p \text{ a sing-}\\ 
\text{singularity of } \ce}} (u_p v_p / c_p) \deg_{D_p} \tilde \ce_r\\
        & = & \sum_{\substack{p \text{ of}\\ \text{ type } \alpha_i}}
 (u(i) v(i) / c_i) +
  \sum_{\substack{p \text{ of }\\ \text{type } \gamma_j}}
 (j (r-j)/ d_j).
\end{eqnarray*}
The second line follows because $\deg_{D_p} \tilde \ce_r = 1$ and
 because over $\gamma_j$ we have $u=j$ and $v=r-j$.  This proves the
 theorem for $<\tilde \ce_r, \mathfrak E_r>$, and the result for $s$ is just the pullback of the relation for $<\tilde 
\ce_s, \mathfrak E_s>$ from $\pic \overline{\stack{S}}^{1/s}_g$. 

\end{proof}

Now we can prove the main theorem.

\begin{proof}(of Theorem \ref{mainthm})
  Since $\mathfrak E_r \otimes \tilde{\ce}^{\otimes r}_r = 
  \omega$, we have 

\begin{eqnarray*}
r<\te,\me> & = & -<\te^{\otimes r}, \me^{-1}> - <\te^{\otimes r}, 
\omega> + <\te^{\otimes r},\omega>\\ & =& <\te^{\otimes r},\omega>-<\te^{\otimes r},\te^{\otimes r}>\\ & 
=& <\te^{\otimes r},\omega>-r^2<\te,\te> 
\end{eqnarray*} 

Deligne-Riemann-Roch and $\mu:=\det f_! \te$ now give
 
\begin{eqnarray*}
&=& <\te^{\otimes r},\omega> -r^2(2\mu-2\lambda+<\ce,\omega>) \\ 
&=& 2r^2 
\lambda -2r^2 \mu-(r-1)<\te^{\otimes r},\omega>\\
&=& 2r^2 \lambda-2r^2 \mu-(r-1)<\omega \otimes \me^{-1},\omega>\\ 
&=& 2r^2 \lambda-2r^2 \mu -(r-1)<\omega,\omega>\\ &=& 
2r^2\lambda-2r^2\mu-(r-1)(12\lambda-\delta)\\ 
\end{eqnarray*} 
Here the last equality follows from the well-known Mumford 
isomorphism: $<\omega,\omega>=12 \lambda-\delta$ (see 
\cite{harris-mumford}). 
\end{proof}

     \subsection{Torsion in $\pic \stacksgr$}\label{tors}

The Picard group of $\stackmg$ is known to be freely generated by 
$\lambda$ when $g$ is greater than $2$ (see 
\cite{arbarello-cornalba:mg-pic}).  And Harer 
\cite{harer:spin-pic-rat} has shown that for $r=2$, the rational 
Picard group $\pic \stack{S}^{1/2}_g \tensor 
\mathbb Q$ has rank one for $g\ge 9$.  So one might expect that $\pic
\stack{S}^{1/2}_g$ is freely generated by $\mu$ or $\lambda$, but Cornalba showed
in \cite{cornalba:spin-pic} that $\pic \stacksgtwo$ has 
$4$-torsion, and one of the consequences of Theorem \ref{mainthm} 
is that whenever $2$ or $3$ divides $r$, there are torsion 
elements in $\pic 
\stacksgr$.  In particular, the following proposition holds.

\begin{proposition}\label{torsion}
If $r$ is not relatively prime to $6$ and $g>1$, then $\pic 
\stacksgr$ has torsion elements: 
\begin{enumerate}
\item If $r$ is even, then $r^2\mu - (r^2 -6r +6)\lambda \neq 0$, and
thus $\frac{1}{2} (r^2\mu - (r^2 -6r +6)\lambda)$ is an element 
of order $4$ in $\pic \stacksgr$. 
\item If $3$ divides $r$, then $\frac23 (r^2\mu - (r^2 -6r +6)\lambda)
\neq 0$, and thus $\frac13 (r^2\mu - (r^2 -6r +6)\lambda)$ is an
element of order $3$ or $6$. 

\item If $r=s d$ and $d$ is even, then $r^2 (\mu_s -\mu_r) - 6 (d^2 -
r d + r + 1) \lambda \neq 0$, and thus $\frac12 r^2 (\mu_s -\mu_r) - 3 (d^2 -
r d + r + 1) \lambda$ has order $4$.
\item If $r=s d$ and $3$ divides $d$, then $\frac23 r^2 (\mu_s -\mu_r) - 4 (d^2 -
r d + r + 1) \lambda \neq 0$, and thus $\frac13 r^2 (\mu_s -\mu_r) - 2 (d^2 -
r d + r + 1) \lambda$ has order $3$ or $6$.
\end{enumerate}
\end{proposition}

\begin{corollary}
 If $6$ divides $r$ then $\frac16(r^2\mu - (r^2 -6r +6)\lambda)$
is an element of order $12$.
\end{corollary}

\begin{proof}[Proof of Proposition~\ref{torsion}]
If the proposition were false and the element in question were 
zero, then in $\pic \stacksgrbar$ this element would be a sum of 
boundary divisors.  In particular, the element in question would 
be of the form $\sum e^{(a,b)}_i \alpha^{(a,b)}_i + \sum 
c^{(a)}_{k,\rho} \gamma^{(a)}_{k,\rho}$.  In the first case, 
multiplication by two, and in the second case, multiplication by 
three, allows us to replace the element in question with a sum 
(from Theorem~\ref{mainthm}) consisting exclusively of boundary 
divisors.  Thus for the first case we have a relation between 
boundary divisors where the sum on the right has all coefficients 
divisible by two: 
$$(1-r) \delta + r <\tilde \ce, \mathfrak E> = 2 \sum e^{(a,b)}_i \alpha^{(a,b)}_i +2 \sum c^{(a)}_{k,\rho} \gamma^{(a)}_{k,\rho}.$$
And for the second case, the
sum on the right has all coefficients divisible by three:
$$(1-r) \delta + r <\tilde \ce, \mathfrak E> = 3 (\text{boundary divisors})$$

Thus by Theorem \ref{lag-indep} the coefficients on the left must also
be divisible by $2$ or $3$, respectively.  However, in both cases 
 the coefficient of $\gamma^{(a)}_{0,\rho}$ for any $\rho$ and $a$ is $1-r$,
which has no divisors in common with $r$; a contradiction.  

The third and fourth cases are similar, but first we must subtract
$d^2$ times the second equation of Theorem~\ref{mainthm}.bis from the
first.  Now reduction mod $2$ and mod $3$ give the necessary
contradictions in the third and fourth cases, respectively. 

\end{proof}

\section{Examples}\label{elliptic}
\subsection{Genus $1$ and Index $1$}

In the case of $g=1$, the only boundary divisors of the stack 
$\overline{\stack{S}}^{1/r}_1:=\overline{\stack{S}}^{1/r, 
\mathbf{0}}_{1,1}$ are those lying over $\delta_0$; that is, the 
Ramond divisors $\gamma_{0,\rho}$, corresponding to the different 
gluings $\{\rho\}$ of $\co_{\mathbb{P}^1}$ at the unique node, 
and the Neveu-Schwarz divisors $\gamma_j$.  In particular, if we 
write $\overline{\stack{S}}^{1/r}_{1}$ as  the disjoint union of 
its irreducible components 
$$\overline{\stack{S}}^{1/r}_{1}= \coprod_{d|r}\overline{\stack{S}}^{1/r,(d)}_{1}$$ where 
a
generic geometric point of $\overline{\stack{S}}^{1/r,(d)}_{1}$ 
has $\ce_r$ isomorphic to  a $d$th root of $\co$, then the only 
boundary divisor in $\overline{\stack{S}}^{1/r,(1)}_{1}$ is 
$\gamma_{0,1}$ corresponding to the trivial bundle $\co_X =\ce_r$ 
on the singular curve $X$.  Moreover, 
$\overline{\stack{S}}^{1/r,(1)}_{1} \rTo 
\overline{\stack{M}}_{1,1}$ is unramified, and so $\delta=\delta_0=\gamma_0=\gamma_{0,1}$
in $\pic \overline{\stack{S}}^{1/r,(1)}_{1}$.  The boundary 
divisor $<\tilde{\ce}, \mathfrak E>$ reduces to zero in this 
case, and so Theorem \ref{mainthm} gives 
\begin{eqnarray*}
2r^2 \mu& =&(2r^2-12r+12)\lambda+(r-1) \delta\\  &= & 
(2r^2-12r+12)\lambda +(r-1)\gamma_0. 
\end{eqnarray*}

We can give a much more complete description of the Picard group 
in this case using Edidin and Graham's equivariant intersection 
theory \cite{edidin-graham} and the following explicit 
construction of $\stack{S}^{1/r,(1)}_{1}$. 

\begin{proposition}\label{quotient}
The stack $\stack{S}^{1/r, (1)}_{1}$ is isomorphic to the 
quotient $\mathbb{A}^2_{c_4,c_6} 
-\Delta:=\{(c_4,c_6)|c^3_4-c^2_6 \neq 0\}$ by the action of $\mathbb{G}_m$, 
given by $v \cdot (c_4,c_6)=(v^{-4r}c_4,v^{-6r}c_6)$.  Similarly,  
$\overline{\stack{S}}^{1/r,(1)}_{1}$ is 
$\mathbb{A}^2_{c_4,c_6}-D:=\{(c_4,c_6)|c_4 \neq 0 \neq c_6\}$ 
modulo the same action of $\mathbb{G}_m$. 
\end{proposition} 

\begin{proof}
First recall that if $\Delta$ is the locus $\{c^3_4-c^2_6\}$ and 
$D$ is the locus $\{c_4=c_6=0\}$, then the stacks 
$\stack{M}_{1,1}$ and $\overline{\stack{M}}_{1,1}$ have a 
representation as the space of cubic forms 
$\{y^2=x^3-27c_4x-54c_6\}$, that is 
$\mathbb{A}^2_{c_4,c_6}-\Delta$ and $\mathbb{A}^2_{c_4,c_6}-D$ 
respectively, modulo the ``standard'' $\mathbb{G}_m$ action $v 
\cdot (c_4,c_6) =(v^{-4}c_4,v^{-6}c_6)$ \cite[Remark following \S 5.4]{edidin-graham}. 
We denote this action by $s:\mathbb{G}_m \cross \mathbb{A}^2 \rTo  
\mathbb{A}^2$ and the action of the proposition by $b: 
\mathbb{G}_m \cross \mathbb{A}^2 \rTo \mathbb{A}^2$.

We have commutative diagrams of stacks
$$\begin{diagram}
(\mathbb{A}^2-\Delta)/b & \rTo & \stack{S}^{1/r,(1)}_{1}\\
\dTo & & \dTo\\
(\mathbb{A}^2-\Delta)/s & \irightarrow & \stack{M}_{1,1}\\
\end{diagram}$$

and

$$
\begin{diagram}
(\mathbb{A}^2 -D)/b & \rTo & 
\overline{\stack{S}}^{1/r,(1)}_{1}\\
\dTo & & \dTo \\
(\mathbb{A}^2-D)/s & \irightarrow & \overline{\stack{M}}_{1,1}\\
\end{diagram}$$
where the top morphism is given by the fact that there is a 
$b$-equivariant choice of a line bundle $\ce_r$ on the curve 
$y^2=x^3-27c_4-54c_6 \subseteq \mathbb{P}^2 \cross (\mathbb{A}^2 
-D)$ and a $b$-equivariant isomorphism $\ce^{\otimes r}_r 
\irightarrow 
\omega$.  The bundle $\ce_r$ is generated by an $r$th root of 
$dx/y$, the invariant differential.  
 This makes sense because $dx/y$ has no zeros or poles. 
  
 Alternately, we may simply take the trivial line bundle $\cn$ on 
 $\mathbb{A}^2$ generated by an element $\zeta$, with the action 
 $b$ defined as $v \cdot \zeta = v^{-1}\zeta$. If $\pi: \{y^2=x^3-27c_4-54c_6\} \rTo \mathbb{A}^2$ is the 
 projection to $\mathbb{A}^2$, then we define 
 $\ce_r$ to be $\pi^* \cn$, and the homomorphism $\ce^{\otimes 
 r}_r \rTo \omega$ to be $\zeta^r \mapsto \frac{dx}{y}$.  This 
 homomorphism is $b$-equivariant  
 since $v \cdot (\frac{dx}{y}) =v^{-r}\frac{dx}{y}$; and it is 
 well-known that for this family $\frac{dx}{y}$ generates the 
 relative dualizing sheaf $\omega_{\pi}$.
 
 The proposition now follows since both the left and right-hand 
 vertical morphisms are \'etale of degree $1/r$, and the bottom 
 morphism is an isomorphism.  Thus the top 
 morphism is \'etale of degree $1$.  
 It is clearly an isomorphism on geometric points, thus also an 
 isomorphism of stacks.
\end{proof}

It is easy to see that the line bundle $\cn$ induces the 
pushforward $\pi_* \ce_r$ on $\overline{\stack{S}}^{1/r,(1)}_{1}$ 
and $\stack{S}^{1/r,(1)}_{1}$.  We denote this bundle by  
$\mu^+$.  Similarly the bundle $\mu^-:=R^1 
\pi_*\ce_r=-\pi_*(\omega 
\otimes \ce_r^{-1})$ can be seen to be $-\lambda+\mu^+$, thus $\mu=\mu^+-\mu^-=\lambda$, 
which is compatible with $12\lambda=\delta$ and with Theorem 
\ref{mainthm}.  Moreover, the explicit description of $\cn$ shows 
that $r \mu^+ 
=\lambda$.  So the order of $\mu^+$ is $12r$ in $\pic 
\stack{S}^{1/r,(1)}_{1}$. 

\begin{corollary}
The Chow rings $A^* (\stack{S}^{1/r,(1)}_1)$ and $A^* 
(\overline{\stack{S}}^{1/r,(1)}_1)$ are isomorphic to 
$\mathbb{Z}[t]/12rt$ and $\mathbb{Z}[t]/24r^2t^2$, respectively.   
Consequently,  $\pic \stack{S}^{1/r,(1)}_{1} = <\mu^+> 
\cong \mathbb{Z}/12r\mathbb{Z}$ and $\pic 
\overline{\stack{S}}^{1/r,(1)}_1 =<\mu^+> \cong \mathbb{Z}$.
\end{corollary}

\begin{proof}
By \cite[Prop. 18 and 19]{edidin-graham} for any smooth quotient 
stack $\cf = [X/G]$ the Chow ring $A^* (\cf)$ is the  equivariant 
Chow ring $A^*_G (X) \cong A^G_* (X)$, and $\pic \cf 
= A^1_G(X)$.  Thus it suffices to compute $A^G_* (X)$, where $X$ is the $\mathbb{A}^2_{c_4,c_6}-\Delta$ or 
$\mathbb{A}^2_{c_4,c_6}-D$, and $G$ in $\mathbb{G}_m$ with the 
action $b: \mathbb{G}_m \cross X \rTo X$ of Proposition 
\ref{quotient}. 

Choosing an $N+1$-dimensional representation $V$ of $G$ with all  
weights $-1$, and letting $U=V-\{0\}$ be the open set where $G$ 
acts freely, then the diagonal action of $G$ on 
$(\mathbb{A}^2_{c_4,c_6}-\{0\}) \cross U$ is free, and $A^G_i 
(\mathbb{A}^2_{c_4,c_6}-\{0\})$ is defined \cite[Defn 
1]{edidin-graham}to be the usual Chow group $A_i 
((\mathbb{A}^2-\{0\})\cross U/G)$ of the quotient scheme 
$((\mathbb{A}^2-\{0\}) \cross U)/G$, which is isomorphic to the 
complement of the zero section of the vector bundle $\co (4r) 
\oplus \co(6r)$ over $\mathbb{P}^N$.  Thus $A^*_G (\mathbb{A}^2_{c_4,c_6}-\{0\}) = \mathbb{Z}[t]/24r^2t^2$.

Moreover, since the form $c^3_4-c^2_6$ has weighted degree $-12r$ 
with respect to the action, the $G$ equivariant fundamental class 
$[\Delta]_G$ of $\Delta$ is $12rt$, and 
$A^G_*(\mathbb{A}^2-\Delta)=A^G_*((\mathbb{A}^2 
-\{0\} \cross U)/[\Delta]_G =\mathbb{Z}[t]/12rt.$ 
Similarly, the class $[D]_G$ is the intersection of $[c_4=0]_G$ 
and $[c_6=0]=(4rt)\cdot(6rt)=0$.  The theorem follows.
\end{proof}

\subsection{Other Components in Genus $1$}

Because we have no explicit representation of 
$\stack{S}^{1/r,(d)}_1$ as a quotient stack for $d>1$, this case 
is more difficult.  Moreover, $\ce_r$ has no global sections, nor 
any higher cohomology, so the bundle $\mu$ is trivial (and 
$\mu^+=\mu^-=0$).  The only other obvious bundles on the stack 
are those induced by pullback along $[d]:\stack{S}^{1/r,(d)}_1 
\rTo \stack{S}^{d/r,(1)}_1$ from $\pic \stack{S}^{d/r,(1)}_1$.
  In particular, we have the bundles $\mu^{d/r,+}$ and $\lambda$.  
  With $\lambda=r\mu^{d/r,+}/d$ and since $12\lambda 
= \delta$, the relation of Theorem \ref{mainthm} gives $$2r^2 
\lambda =0 
\text{ in } \pic \stack{S}^{1/r,(d)}_1.$$ In particular, if $6$ 
does not  divide $r$, then the homomorphism $[d]^*:\pic 
\stack{S}^{d/r,(1)}_{1} \rTo 
\pic \stack{S}^{1/r,(d)}_1$ is not injective.

In the case of $d=2$ or $3$ we can follow Mumford 
\cite{mumford:picard} and construct $r$-spin curves of index $d$ 
which have non-trivial automorphisms,  and these will give 
homomorphisms $\pic 
\stack{S}^{1/r,(d)}_1 
\rTo \mathbb{G}_m$.  In particular, in the case that $d=2$, consider the curve
$E_{1728}:y^2 =x(x^2-1)$, and the two-torsion point $p=(0,0)$.  
Associated to $p$ is the line bundle $\ce_r:=\{f\psi|f\in 
k(E_{1728}), 
\, (f) \geq p- \infty\}$, and the isomorphism $c_{r,1}:\ce_r^{\otimes r} \irightarrow \omega$ 
defined by $$g \psi^r \mapsto \frac{gdx}{x^{r/2}y}.$$ This is an 
isomorphism because $x$ is a global section of $\co(-2p+2\infty)$ 
giving an isomorphism to $\co$, $\frac{dx}{y}$ is a  global 
section of $\omega$ giving an isomorphism to $\co$, and $2$ 
divides $r$. 

An automorphism $\sigma$ of order $4$ of the underlying curve 
$E_{1728}$ 
 can be defined as $\sigma(x,y)=(-x,iy)$, and 
$\sigma$ can be extended to an automorphism of the entire spin 
curve $(E_{1728},(\ce_r,c_{r,1}))$ by letting $\zeta$ be a 
primitive $4r$th root of unity such that $\zeta^r =(-1)^{r/2}i$ 
and then defining 
$$\sigma(\psi) =\zeta \psi.$$ $c_{r,1}$ is preserved by $\sigma$ 
since $\sigma (\frac{dx}{y x^{r/2}})=(-1)^{r/2}i(\frac{dx}{y 
x^{r/2}}) = \zeta^r (\frac{dx}{y x^{r/2}}).$

$(\ce_r,c_{r,1})$ is clearly of index $(2)$, and although $\pi_* 
\ce_r ,\, R^1 \pi_*\ce_r$, and $\pi_! \ce_r$ are all trivial (zero or $\co$),
on $\stack{S}^{1/r,(2)}_1$, the bundle $\ce_{r/2}=\ce_r$ is 
always isomorphic to $\co_E$, and so the sheaves 
$\mu^{2/r,+}:=\pi_*\ce_{r/2}$ and $ 
\mu^{2/r,-}:=R^1\pi_*\ce_{r/2}$  are always line bundles on 
$\stack{S}^{1/r,(2)}_1$ (actually they are just the pullbacks of 
$\mu^+$ and $\mu^-$, respectively, from $\stack{S}^{2/r,(1)}_1$ 
via the morphism $[2]:\stack{S}^{1/r}_1 
\rTo \stack{S}^{2/r}_1$).

For every line bundle $\cl \in \pic \stack{S}^{1/r,(2)}_1$, the 
geometric point $(E_{1728},(\ce_r, c_{r,1}))$ of 
$\stack{S}^{1/r,(2)}_1$  associates a one-dimensional vector 
space $\cl(E_{1728},(\ce_r,c_{r,1})) 
\cong k$, and $\sigma$ induces an automorphism $\cl (\sigma) \in 
\mathbb{G}_m$ of $\cl(E_{1728},(\ce_r,c_{r,1}))$.  Moreover, 
$\cl(\sigma)$ has order dividing $4r$ (the order of $\sigma$).  
Thus we have a homomorphism $$ \pic \stack{S}^{1/r,(2)}_1 \rTo 
<\zeta> \cong \mathbb{Z}/4r.$$ Moreover, the explicit description  
of $\ce_r$ shows that $\mu^{2/r,+}$ maps to $\zeta^2$ and 
$\mu^{2/r,-}$ maps to $\zeta^{2-r}$. 

In the case of $d=3$ we can use a similar argument, applied to 
the curve $E_0:y^2+y=x^3$, the line bundle $\ce_r=\{f \psi|f \in 
k(E_0),\,(f) \geq p - \infty\}$ and the isomorphism $$ c_{r,1} 
:\psi^r \mapsto \frac{dy}{x \cdot y^{r/3}}.$$  Choose a primitive $3r$th root of unity 
$\xi$,
and define the automorphism $$\tau:(x,y) \mapsto (\xi^{-r}x,y)$$ 
$$\psi \mapsto \xi \psi.$$
This is compatible with $c_{r,1}$, and so is an automorphism of 
$(E_0,(\ce_r,c_{r,1}))$.  This gives a homomorphism $\pic 
\stack{S}^{1/r,(3)}_1 \rTo <\xi> \cong \mathbb{Z}/3r$ and the elements $\mu^{3/r,+}:= \pi_*
\ce_{r/3} (=[3]^*\mu^+)$ and $\mu^{3/r,-}:=R^1 \pi_* \ce_{r/3}$ map to $\xi^3$
 and $\xi^{3-r}$, respectively. Thus we have the following 
 commutative diagrams.
 
 \begin{diagram}
 \pic \stack{S}^{3/r,(1)}_1 & \irightarrow & 
 <\mu^{3/r,+}>=\mathbb{Z}/12r\\
 \dTo & & \dTo\\
 \pic \stack{S}^{1/r,(3)}_1 &\rTo& \mathbb{Z}/3r 
 \end{diagram}
 
 and
 
 \begin{diagram}
 \pic \stack{S}^{1/r,(1)}_1 & \irightarrow & 
 <\mu^{2/r,+}>=\mathbb{Z}/12r\\
 \dTo & & \dTo\\
 \pic \stack{S}^{1/r,(2)}_1 & \rTo & \mathbb{Z}/4r
 \end{diagram}
where the right-hand vertical morphisms take $\mu^{3/r,+}$ to $3$ 
and $\mu^{2/r,+}$ to $2$, respectively, and $\lambda = r 
\mu^{3/r,+}$ or $r\mu^{2/r,+}$, respectively.

For $d>3$ this strategy does not work as  well, since no 
automorphisms of the underlying curve preserve a spin structure 
of type $d$.  Nevertheless, we still have for any 
$(E,(\ce_r,c_{r,1}))$ the automorphism defined by taking $\ce_r$ 
to $\eta \ce_r$, where $\eta$ is any $r$th root of unity.  This 
gives, for $\eta$ a primitive $r$th root of unity, $$\pic 
\stack{S}^{1/r,(d)}_1 \rTo <\eta> = \mathbb{Z}/r$$ and the bundle 
$\mu^{d/r,+}:=\pi_* \ce_{r/d}$ maps to $\eta^d$. 

This inspires the following conjecture:

\begin{conjecture}
$$\pic \stack{S}^{1/r,(d)}_1 = <\mu^{d/r,+}> = \begin{cases}
\mathbb{Z}/2r & \text{ if } d=2 \\
\mathbb{Z}/r  & \text{ if } d=3\\
\mathbb{Z}/(r/d) & \text{ if } d>3
\end{cases}$$
\end{conjecture} 
 
\subsection{General genus, $r=2$}
For $g>2$, and $r=2$, we have that $2 \mu + \lambda$ is an 
element of order $4$, and $\lambda$ (and hence $\mu$) is an 
element of infinite order.  Moreover, Harer has proved  for $g>9$ 
that $H_1(\stack{S}^{1/2,even}_g,\mathbb{Z}) 
= H_1 (\stack{S}^{1/2,odd}_g,\mathbb{Z}) = \mathbb{Z}/4$, and 
$H^2(\stack{S}^{1/2,even}_g,\mathbb{Q}) = H^2 
(\stack{S}^{1/2,odd}_g, \mathbb{Q})= \mathbb{Q}$, so for $g>9$ we  
have $$\pic \stack{S}^{1/2,odd}_g  = \pic \stack{S}^{1/2,even}_g 
\cong 
\mathbb{Z} \cross 
\mathbb{Z}/4\mathbb{Z}.$$

What is not yet completely clear, but seems reasonable to expect, 
is the following presentation for that group. 

\begin{conjecture}
For all $g>2$, $$\pic \stack{S}^{1/2,even}_g =\pic 
\stack{S}^{1/2,odd}_g= <\mu,\lambda 
\,|\, 8 
\mu+4 
\lambda =0>.$$
\end{conjecture}

\section*{Conclusion}

We have worked out many  relations between the elements of $\pic 
\stacksgrbar$ and $\pic \stack{S}^{1/r}_g.$  This generalizes the work of Cornalba 
\cite{cornalba:theta,cornalba:spin-pic}, whose results hold  in 
the case where $r=2$.  One of the interesting consequences of 
these relations is the existence of elements in $\pic \stacksgr$ 
of elements of order $4$ if $2$ divides $r$, and elements of 
order $3$, if $3$ divides $r$.  Somehow $2$ and $3$ seem to be 
special, however, and when $r$ is relatively prime to $6$ and 
$g>2$, there do not appear to be any torsion elements in $\pic 
\stacksgr$. 

Corresponding results for $\pic \stacksgrnmbar$ will appear in 
\cite{jkv}, where they are used to prove the genus-zero case of 
the generalized Witten conjecture \cite{witten:r-spin-conj}. 
 
\section*{Acknowledgments}

I wish to thank Takashi Kimura, J\'anos Koll\'ar, Bill Lang, and 
Arkady Vaintrob for helpful discussions and suggestions regarding 
this work.  I am also grateful to Heidi Jarvis for help with 
typesetting. 

\providecommand{\bysame}{\leavevmode\hbox to3em{\hrulefill}\thinspace}

\end{document}